\documentclass[reqno, 12pt]{amsart}
\usepackage{amsfonts}
\usepackage{amssymb,amsmath}
\usepackage{amsthm}
\usepackage{upref}

\makeatletter
\def\LaTeX{\leavevmode L\raise.42ex
    \hbox{\kern-.3em\size{\sf@size}{0pt}\selectfont A}\kern-.15em\TeX}
 \makeatother

\numberwithin{equation}{section}

\newtheorem{lemma}{Lemma}[section]
\newtheorem{theorem}[lemma]{Theorem} 
\newtheorem{corollary}[lemma]{Corollary}
\newtheorem{proposition}[lemma]{Proposition}

\theoremstyle{definition}

\newtheorem{definition}[lemma]{Definition}
\newtheorem{example}[lemma]{Example}

\newtheorem{assumption}[lemma]{Assumption}

\newtheorem{remark}[lemma]{Remark}

 \newcommand{\supp}{\operatorname{supp}}
  \newcommand{\e}{\eqref}

\newcommand{\q}{\quad}

\newcommand{\ti}{\tilde}
\newcommand{\wt}{\widetilde}

\newcommand{\la}{\langle}
\newcommand{\ra}{\rangle}

\newcommand{\ov}{\overline}
 \renewcommand{\d}{\delta}

\renewcommand\Im{\operatorname{Im}}
\renewcommand\Re{\operatorname{Re}}

\newenvironment{pf}{\begin{proof}}{\end{proof}}

\def\qqq{\mathrel{\subset\mkern-15mu\lower.38ex\hbox{${\scriptscriptstyle\rightarrow}$}}}

\let\cal\mathcal

\let\Bbb\mathbb

    \DeclareMathOperator{\Res}{Res}
  
\begin{document}

\title   {Quasi-diagonalization of Hankel  operators}
\author{ D. R. Yafaev}
\address{ IRMAR, Universit\'{e} de Rennes I\\ Campus de
  Beaulieu, 35042 Rennes Cedex, FRANCE}
\email{yafaev@univ-rennes1.fr}


\begin{abstract}
We show that all Hankel operators $H$ realized as integral operators with kernels $h(t+s)$ in $L^2 ({\Bbb R}_{+}) $  can be quasi-diagonalized as $H= {\sf L}^* \Sigma {\sf L} $. Here ${\sf L}$ is
the Laplace transform,  $\Sigma$ is the operator of multiplication by a function (distribution) $\sigma(\lambda)$, $\lambda\in {\Bbb R}$. 
We find a scale of spaces of test functions  where ${\sf L} $ acts as an isomorphism. Then ${\sf L}^*$ is an isomorphism of the corresponding spaces of distributions. We show that $h= {\sf L}^* \sigma$ which yields a one-to-one correspondence between kernels $h(t)$  and  sigma-functions $\sigma(\lambda)$ of Hankel operators. The sigma-function   of a self-adjoint Hankel operator $H$ contains  substantial information about its spectral properties. Thus we show that the operators 
$H$ and $\Sigma$ have the same numbers of positive and negatives eigenvalues. In particular, we find necessary and sufficient conditions for sign-definiteness of Hankel operators. These results are illustrated at   examples of quasi-Carleman operators
generalizing the classical Carleman operator with kernel $h(t)=t^{-1}$ in various directions. The concept of the sigma-function directly leads to a criterion  (equivalent of course to the classical Nehari theorem) for boundedness 
of Hankel operators. 
Our construction also shows that every Hankel operator is unitarily equivalent by the Mellin transform to a pseudo-differential operator with amplitude which is a product of functions of one variable only (of $x\in{\Bbb R}$ and of its dual variable). 
  \end{abstract}

\maketitle


\section{Introduction}  

{\bf 1.1.}
Hankel operators can be defined as integral operators  
\begin{equation}
(H f)(t) = \int_{0}^\infty h(t+s) f(s)ds 
\label{eq:H1}\end{equation}
in the space $L^2 ({\Bbb R}_{+}) $ with kernels $h$ that depend  on the sum of variables only. We refer to the books \cite{NK,Pe,Po} for basic information on Hankel operators. Of course $H$ is symmetric if $  h(t)=\ov{h(t)}$. There are very few cases when Hankel operators can be explicitly diagonalized. The   simplest and most important case $h(t)=t^{-1}$ was considered by T.~Carleman in \cite{Ca}.

Our goal here is to show that all Hankel operators can be quasi-diagonalized in the following sense. Let ${\sf L}$,
   \begin{equation}
  ({\sf L} f) (\lambda)= \int_{0}^\infty e^{-t \lambda} f(t) dt  ,
\label{eq:LAPj}\end{equation}  
be the Laplace transform. Under very general assumptions on  $h$, we prove that
\begin{equation} 
H= {\sf L}^* \Sigma {\sf L} 
 \label{eq:MAIDp}\end{equation} 
 where $\Sigma$ is the operator of multiplication by the function $\sigma(\lambda)$ {\it formally}  linked to the kernel $h$ by the relation
    \begin{equation}
  h (t)  =    \int_{-\infty}^\infty     e^{-t \lambda} \sigma(\lambda)  d \lambda
   \label{eq:conv1}\end{equation}
   (that is, $h$ is the two-sided Laplace transform of $\sigma$).  We call $\sigma(\lambda)$ the {\it sigma-function} of the Hankel operator $H$ or of its kernel $h(t)$.

It is clear from   formula \e{eq:conv1} that $\sigma(\lambda)$  can be a regular function only for   kernels $h(t)$      satisfying  some specific analytic  assumptions. Without such very restrictive assumptions,  $\sigma$ is necessarily a distribution. Even for very good kernels $h(t)$ (and especially for them), $\sigma(\lambda)$ may be a highly singular distribution.
For example,    for $h(t)=t^k e^{-\alpha t}$ where  $\Re \alpha >0$ ($\alpha$ may be complex) and $k=0,1, \ldots$, the sigma-function $\sigma(\lambda)= \d^{(k)}(\lambda-\alpha )$ is a  derivative   of the delta-function.

Relation \e{eq:MAIDp} does not require the condition $  h(t)=\ov{h(t)}$. If however it is satisfied, then under proper assumptions $H$ can be realized as a self-adjoint operator although $\Sigma$ is determined by a quadratic  form which does not necessarily give rise to a self-adjoint operator.

   Let us compare quasi-diagonalization \e{eq:MAIDp} of Hankel operators with the standard diagonalization of convolution operators $B$ with integral kernels $b(x-y)$ in the space $L^2 ({\Bbb R}) $. Let $\Phi$ be the Fourier transform, and let $S$ be the operator of multiplication by the function (the symbol of the convolution operator $B$) $s (\xi)= \sqrt{2\pi} (\Phi b)(\xi) $, $\xi\in {\Bbb R}$. Then 
   \begin{equation} 
B=\Phi^* S \Phi .
 \label{eq:MAIDq}\end{equation} 
  Since the operator $\Phi$ is unitary, formula 
\e{eq:MAIDq} reduces  convolution operators to      multiplication operators and hence exhibits their complete spectral analysis. 

This is not of course the case with Hankel operators because ${\sf L} $ is not unitary. Fortunately, in an appropriate sense,  ${\sf L} $ turns out to be invertible. Therefore it follows from relation \e{eq:MAIDp} that, in the self-adjoint case,  the total numbers of (strictly) positive $N_{+} (H) $ and negative $N_- (H) $ eigenvalues of a Hankel operator $H$ equal the same quantities for the operator $\Sigma$ of multiplication by the  function $\sigma(\lambda)$:
  \begin{equation}
N_{\pm}(H)= N_{\pm}(\Sigma).
\label{eq:MAIDT}\end{equation} 
In particular,      $\pm H\geq 0$ if and only if $\pm \Sigma\geq 0$. Moreover, if $\sigma (\lambda) >0$ (or $\sigma (\lambda) <0$) on a set of positive Lebesgue measure, then the Hankel operator $H$ has infinite positive (or negative) spectrum. On the other hand, singularities of $\sigma(\lambda)$ at some isolated points produce   finite numbers (depending on the order of the singularity) of positive or negative eigenvalues.

  Equality  \e{eq:MAIDT} can be compared with Sylvester's inertia theorem  which states the same   for Hermitian matrices $H$ and $\Sigma$ related by equation \e{eq:MAIDp} provided the matrix $\sf L$ is invertible. In contrast to the linear algebra, in our case the operators $H$ and $\Sigma$ are of a completely different nature and $\Sigma$ (but not $H$) admits an explicit spectral analysis.

  Hankel operators can also be realized in the space     $ l^2({\Bbb Z}_{+})$ of sequences $g =(g_{0}, g_{1}, \ldots)$ 
   by the relation 
\begin{equation}
(Q g)_{n} =   \sum_{m=0}^\infty q_{n+m}g_{m}, \q g=(g_{0}, g_{1}, \ldots),
\label{eq:K1}\end{equation}
which is obviously a discrete analogue of continuous definition \e{eq:H1}. So it is not astonishing that there exists a unitary operator ${\bf U} : l^2({\Bbb Z}_{+})\to  L^2({\Bbb R}_{+})$ such that the operator
\begin{equation}
H={\bf U} Q{\bf U} ^{-1}
\label{eq:K2}\end{equation}
   acting in $L^2({\Bbb R}_{+})$ is Hankel if and only if $Q$ is a Hankel operator in $l^2({\Bbb Z}_{+})$. However the
 construction of  
   the operator ${\bf U} $ is nontrivial and is given in terms of the Laguerre functions. 
   
   Our results can be translated into the space  $ l^2({\Bbb Z}_{+})$. In particular, it follows from \e{eq:conv1} 
   that
         \begin{equation}
 q_{n } =   \int_{-1}^1 \eta(\mu) \mu^n d\mu
\label{eq:KL5}\end{equation}
where the function $\eta(\mu)$ is linked to the sigma-function $\sigma(\lambda)$ by a simple change of variables. Equations \e{eq:KL5}  for $\eta(\mu)$ are known as the Hausdorff moment problem.
Thus the construction of the sigma-function provides an efficient procedure for the solution of this problem.

 \medskip
 
 {\bf 1.2.} 
 The precise sense of formula \e{eq:MAIDp} needs of course  to be clarified. Actually, instead of \e{eq:MAIDp} we prove the identity    
  \begin{equation}
 (H f_{1}, f_{2}) =  ( \Sigma {\sf L} f_{1}, {\sf L} f_{2})
\label{eq:MAIDs}\end{equation} 
on a suitable space of  test functions $f_{1}, f_{2} $.     We find a scale of   spaces  of  test functions where ${\sf L}$ acts as an isomorphism. By duality, the adjoint operator ${\sf L}^*$  establishes an isomorphism of the corresponding spaces of distributions. Relation \e{eq:conv1} should also be understood in the sense of distributions and, strictly speaking, it means that $h ={\sf L}^* \sigma$, that is, 
 \begin{equation}
  \sigma= ({\sf L}^*)^{-1}h .
  \label{eq:ss}\end{equation} 
  Therefore, instead of operators, we consequently work with quadratic forms  which is both more general and more convenient.
It is natural to  also  treat $h$   as a distribution. This yields a one-to-one correspondence between kernels $h$ of Hankel operators and their sigma-functions $\sigma$ and makes the theory   self-consistent. 

To be precise, equality    \e{eq:MAIDT} is also formulated in terms of the corresponding quadratic forms $ (H f , f ) $ and $( \Sigma w, w)$. If the form $ (H f , f ) $ gives rise to the self-adjoint operator $H$, then  \e{eq:MAIDT} yields an explicit expression for $N_{\pm} (H)$. We emphasize that typically $\Sigma$ cannot be realized as a self-adjoint operator.

Although formula \e{eq:MAIDp} does not give the diagonalization of a Hankel operator $H$, it shows that $H$ is unitarily equivalent (the corresponding unitary transformation is essentially the Mellin transform) to  a pseudo-differential operator  $A$ in the space $L^2 ({\Bbb R})$ with  the amplitude  
  \begin{equation}
  v(\xi) s(x) v(\eta), \q x, \xi, \eta\in {\Bbb R},
 \label{eq:sam}\end{equation}
that  factorizes into a product of functions depending on one variable only.
Here 
   \[
v(\xi) = \frac{\sqrt{\pi} }{\sqrt{\cosh(\pi\xi)} }
\]
 is quite   explicit and $s(x)$ (called the sign-function of a Hankel operator $H$ in Ê\cite{Y})
is linked to the sigma-function by the formula
   \begin{equation}
s(x)=\sigma (e^{-x}).
 \label{eq:sig2}\end{equation}
 To put it  differently, $A$ is the integral operator   in the space $L^2 ({\Bbb R})$ with kernel
   \[ 
(2\pi)^{-1/2}  v(\xi) \hat{s}(\xi-\eta) v(\eta) 
 \]
 where $\hat{s}=\Phi  s$ is the Fourier transform of the sign-function $s$.

 \medskip
 
 {\bf 1.3.} 
   We emphasize that the sigma-function $\sigma (\lambda)$ of a Hankel operator $H$  and its symbol $\theta (\xi)$, $\xi\in {\Bbb R}$,  
    are   different objects. In some sense they are dual to each other.  Let us discuss their link at  a {\it formal} level
    for bounded Hankel operators $H$ when $\sigma (\lambda)=0$ for $\lambda<0$.   It is more convenient for us to work with symbols $\omega(\mu): = \theta  (i\mu)$   defined on the imaginary axis $\Re \mu=0$. 
         Recall that the kernel $h(t)$ and the symbol  $\omega(\mu)$ of a Hankel operator  are related by the formula
    \begin{equation}
    h(t)= (2\pi i)^{-1} \int_{-i\infty}^{i\infty } e^{ t \mu} \omega (\mu) d\mu, \q t>0.
\label{eq:hs}\end{equation}
       Equality \e{eq:hs} does not determine $ \omega (\mu)$ uniquely, but it is satisfied if 
  \begin{equation}
   \omega (\mu)=\int_{0}^\infty e^{-t \mu} h(t) dt.
\label{eq:hst}\end{equation}
   This function is analytic  in the right half-plane. Substituting here \e{eq:conv1}, we see that $\sigma$ and $\omega$ are linked by the Stieltjes transform:
   \begin{equation}
     \omega (\mu)=\int_{0}^\infty  (\lambda+\mu)^{-1} \sigma (\lambda) d\lambda.
    \label{eq:hs1}\end{equation}
    
    The relation between  $\sigma$ and $\omega$ can also be expressed in the following way. Let ${\Bbb H}^1_{r}$ be the Hardy class  of functions analytic in the right half-plane, and let $g \in {\Bbb H}^1_{r}$.   Since (see, e.g., the book \cite{Koos}, page~156)
       \[
     g(\lambda)=\frac{1}{2\pi i}\int_{-i\infty}^{i\infty} \frac{g(\mu)}{\lambda-\mu}d\mu, \q \Re\lambda>0,
     \]
      it   follows from \e{eq:hs1} that
      \begin{equation}
 2\pi i  \int_{0}^{\infty} \sigma(\lambda)  \ov{g(\lambda)}d\lambda  = \int_{-i\infty}^{i\infty} \omega(\mu) \ov{g(\mu)}d\mu .
    \label{eq:hs2a}\end{equation}  
     
     Formula \e{eq:MAIDp}  is of course consistent with the standard representation of a Hankel operator $H$ in terms of its symbol $\omega$.  Indeed, by the Paley-Wiener theorem, the operator $(2\pi)^{-1/2}{\sf L}$ is a unitary mapping of $L^2 ({\Bbb R}_{+}) $ onto the Hardy class ${\Bbb H}^2_{r}$. For $f \in L^2 ({\Bbb R}_{+}) $ and $\Re \mu\geq 0$,  we put $\ti{f}( \mu)=(2\pi)^{-1/2} ({\sf L}f) ( \mu)$. Then, by the definition of the symbol, we have
       \begin{equation}
  (Hf,f)      = -i  \int_{-i\infty}^{i\infty} \omega(\mu) 
     \ti{f} (-\mu) \ov{\ti{f} ( \mu)}d\mu.
    \label{eq:hs2b}\end{equation} 
    Let us now apply relation \e{eq:hs2a} to the  function 
    $g (\mu)= \ov{\ti{f} (\bar{\mu})} \ti{f} ( \mu)$. Then putting together formulas \e{eq:hs2b}  and 
          \[
   ( {\sf L}^* \Sigma {\sf L} f , f )     = 2\pi \int_{0}^\infty \sigma(\lambda)
    | \ti{f} (\lambda)|^2 d\lambda,
    \]
     we recover  representation \e{eq:MAIDp}.

  Notice   that, in contrast to the symbol, the sigma-function contains substantial information about spectral properties of $H$. Comparing relations \e{eq:MAIDp} and \e{eq:MAIDq}, we argue that in the theory of Hankel operators, it is rather sigma-functions (and not symbols) that play the role of symbols of convolution operators. We also note that there is the one-to-one correspondence between kernels and sigma-functions and that the notion of the sigma-function does not require the boundedness of $H$.

 \medskip
 
 {\bf 1.4.}
 A {\it formal} proof   of  the identity \e{eq:MAIDp}
 is quite simple and is actually the same as that of the identity \e{eq:MAIDq} for convolutions. Indeed, the integral kernel of the operator in the right-hand side of  \e{eq:MAIDp} equals
 \[
 \int_{-\infty}^\infty e^{-\lambda t} \sigma (\lambda) e^{-\lambda s} d \lambda=h(t+s)
 \]
 if $ \sigma (\lambda)$ and $h(t)$ are linked by formula  \e{eq:conv1}. Thus it equals the integral kernel of the Hankel operator  $H$.
 
  However a rigorous proof of \e{eq:MAIDp} or, more precisely, of \e{eq:MAIDs} requires a choice of a suitable set of test functions $f_{j}(t)$, $j=1,2$,  and a correct formulation of relation  \e{eq:conv1}. The most natural and general choice is to work on functions $f_{j} \in C_{0}^\infty ({\Bbb R}_{+})$ and to require that $h \in C_{0}^\infty ({\Bbb R}_{+})'$. Note that 
  ${\sf L}: C_{0}^\infty ({\Bbb R}_{+}) \to {\cal Y}$ where the set ${\cal Y}$ consists of analytic functions $g(\lambda)$ exponentially decaying as $\Re \lambda \to +\infty$, exponentially bounded as $\Re \lambda \to -\infty$ and decaying faster than any power of
  $|\lambda|^{-1}$ as $|\Im \lambda| \to \infty$. Since   ${\sf L}: C_{0}^\infty ({\Bbb R}_{+}) \to {\cal Y}$ and hence ${\sf L}^*: {\cal Y}' \to C_{0}^\infty    ({\Bbb R}_{+})'$ are isomorphisms, we see that  $\sigma\in {\cal Y}'$.
  
 It turns out that typically regular kernels (like those of finite rank Hankel operators) yield singular sign-functions. On the contrary, singular kernels (such as $h(t)=t^{-q}$ where $q>0$ may be arbitrary large) yield smooth  sign-functions.  Nevertheless the conditions $h \in C_{0}^\infty ({\Bbb R}_{+})'$  and  $\sigma\in {\cal Y}'$ are equivalent and $h$ can be recovered from $\sigma$ by formula \e{eq:conv1}. Thus, although singularities of $h$ and $\sigma$ may be quite different, there is the one-to-one correspondence  between $h$ and $\sigma$ in the  classes of distributions
 $  C_{0}^\infty ({\Bbb R}_{+})'$  and  $  {\cal Y}'$, respectively.
 
 Another possibility is to work on  a set of test functions $f(t)$ satisfying certain analyticity assumptions. This approach is more  symmetric because for such $f$, functions $({\sf L}f)(\lambda)$ satisfy  conditions similar to those on $f(t)$. This leads to the one-to-one correspondence between kernels $h$ and sigma-functions $\sigma$ in the   dual spaces of distributions. It is noteworthy that the inversion of the Laplace transform ${\sf L}$ in the spaces of analytic test functions is quite explicit and
  relies on its factorization.

  In specific examples, the consideration of the form $(\Sigma w, w)$ on analytic functions $w= {\sf L}f\in {\cal Y}$ is not always convenient. Fortunately under mild additional assumptions on the sigma-function $\sigma$,  the set $ {\cal Y}$ of test functions $w$ can be replaced by functions $w\in C_{0}^\infty ({\Bbb R}_{+})$. The proof of this reduction   also relies on the factorization of the Laplace transform ${\sf L}$.

As was already mentioned, even for very regular kernels $h$, the sigma-function 
 $  \sigma$
   may be a highly singular distribution. However,   we show that, for positive\footnote{We always use the term ``positive" (``negative") for a nonnegative (nonpositive) operator or a function. Otherwise we write ``strictly positive" (``negative") } Hankel operators, $\sigma(\lambda)d\lambda$ is given by some positive measure. Thus in the sign-definite case,  $\sigma(\lambda)$ cannot be more singular than delta-functions $\d(\lambda-\alpha)$ where $\alpha>0$. Note that, for {\it positive}  Hankel operators, the concept of the sigma-function, or rather of the associated measure $\sigma(\lambda)d\lambda$, goes back at least to Hamburger (see his paper \cite{Hamb} on moment problems or Theorem~2.1.1 in \cite{AKH}) and to Bernstein (see his theorem on exponentially convex functions in \cite{Bern} or Theorem~5.5.4 in \cite{AKH}). Thus, to a certain extent, our results can be considered as an extension of these classical theorems to the non-sign-definite case.

 \medskip
 
 {\bf 1.5.} 
 We illustrate our general results   at the example of  
  kernels  
\begin{equation}
h(t)= (t+r)^k e^{-\alpha t}, \q r\geq 0, 
\label{eq:E1r}\end{equation}
where $\alpha$ and $k$ are  arbitrary real numbers.   These  kernels give rise to Hankel operators if $\alpha>0$ or $\alpha=0$, $k<0$.  
If   $ \alpha=r=0$ and $k=-1$, then  $H$ is the Carleman operator. In the general case we use the term ``quasi-Carleman operator" for a Hankel operator with kernel   \e{eq:E1r}.
 

  We show that for kernels \e{eq:E1r} the sigma-function  defined by  relation \e{eq:ss} is given by the explicit formula
    \begin{equation}
\sigma (\lambda)=
   \frac{1}{\Gamma(-k)}   (\lambda -\alpha)_{+}^{-k-1}e^{-r (\lambda-\alpha)} , \q k\not\in {\Bbb Z}_{+},
\label{eq:bbr7}\end{equation}
where $\Gamma(\cdot)$ is the gamma function and  $\mu_{+}^{-k-1}$ is the standard distribution defined below by formula \e{eq:di}. If $k\in {\Bbb Z}_{+}$, then $\sigma (\lambda)$ is expressed in terms of the derivatives of the Dirac function:
 \begin{equation}
\sigma (\lambda)= \d^{(k)}(\lambda-\alpha) e^{-r(\lambda-\alpha)}.
\label{eq:dii}\end{equation}
Distributions \e{eq:bbr7} and \e{eq:dii} may be singular at the point $\lambda=\alpha$, and the order of the singularity is determined by the parameter $k$. We show that  the  numbers $N_{\pm} (H)$ are also determined by the parameter $k$ only.

 If $k<0$, then $\sigma \in L^1_{\rm loc} ({\Bbb R}_{+})$ and $\sigma (\lambda)\geq 0$ so that $H\geq 0$.

On the contrary, if $k>0$, then  function \e{eq:bbr7} has the singularity at the point $\lambda=\alpha$ which gets stronger as $k$ increases. If $k\not\in {\Bbb Z}_{+}$, then the function $\sigma(\lambda)$ for $\lambda\neq\alpha$ has the same sign as $\Gamma(-k)$. Therefore $H$ has infinite positive (negative) spectrum if the integer part $[k]$ of $k$ is odd (even). The analysis of the singularity of the function $\sigma(\lambda)$ at the point $\lambda=\alpha$ shows that
 $N_+ (H )=  [k] /2+1$   for even $[k]$ and  $N_- (H) = ([k]+ 1)/2$  for odd $[k]$. 
  If $k\in {\Bbb Z}_{+}$, then the operator $H$ has finite rank $k+1$. In this case it follows from formula \e{eq:dii} (see 
  \cite{Yf}, for details)  that     $N_{+} (H)= N_{-} (H)+1= k/2
+1$ if $k$ is even and
  $N_\pm (H) =(k+1)/2$ if $k$ is odd.
 

We emphasize that, for example,  for $\alpha> 0$, $k>- 1$ and arbitrary $r\geq 0$,  Hankel operators $H$ are  compact, but  $\Sigma$ are not    defined as bounded  operators because of the singularity of the function $\sigma (\lambda)$ at the  point $\lambda=\alpha$.

     \medskip
 
 {\bf 1.6.}
 Let us briefly describe the structure of the paper. Section~2 plays the central role. Here we give the precise definition of the sigma-function, prove the main identity and discuss its consequences. In Section~2, we work on the space $C_{0}^\infty ({\Bbb R}_{+})$ of test functions $f(t)$.   Section~3 is specially devoted to bounded Hankel operators. Here we elucidate the relation between symbols and sigma-functions and prove an analogue of the Nehari theorem in terms of sigma-functions.  We collect various results relying on the factorization of the Laplace transform in Section~4. In particular, we check here that every Hankel operator $H$ is unitarily equivalent to a pseudo-differential operator with amplitude \e{eq:sam}. Then we show that, by a study of the form $(\Sigma w,w)$, a set of analytic test functions $w$ can be replaced by the set $C_{0}^\infty ({\Bbb R}_{+})$. This is technically essentially more convenient. Finally, we     carry over here the results of Section~2 to spaces of analytic test functions $f(t)$.    The case of positive Hankel operators when $\sigma(\lambda)$ is  determined by a measure is discussed in Section~5.  
    Hankel operators $H$ with kernels  \e{eq:E1r} and its various generalizations  are studied in Section~6 where we find an explicit formula for the numbers $N_{\pm}(H)$.  Finally, in Section~7 we   discuss a translation of our results into the representation of Hankel operators in the space $l^2 ({\Bbb Z}_{+})$ of sequences.
   
   
  Let us introduce some standard
 notation.     We denote by $\Phi$, 
\[
(\Phi u) (\xi)=  (2\pi)^{-1/2} \int_{-\infty}^\infty u(x) e^{ -i x \xi} dx,
\]
  the Fourier transform and recall    that  
  $\Phi $ is the one-to-one mapping of the Schwartz  space ${\cal S}={\cal S} ({\Bbb R}) $  onto itself. Moreover, $\Phi$ as well as its inverse   $\Phi^{-1}$ are continuous mappings. In such cases we say that a mapping  is an isomorphism.
   The dual class of distributions (continuous antilinear functionals on ${\cal S}$) is denoted ${\cal S}'$. 
  We use the notation   ${\pmb\la} \cdot, \cdot {\pmb\ra}$ and $\la \cdot, \cdot\ra$ for   the  
  duality symbols in $L^2 ({\Bbb R}_{+})$ and $L^2 ({\Bbb R})$, respectively. They are   linear in the first argument and antilinear in the second argument.

  We often use the same notation for a function and for the operator of multiplication by this function.    
   The letter $C$ (sometimes with indices) denotes various positive constants whose precise values are inessential; $\d_{n,m}$ is the Kronecker symbol, i.e., $\d_{n,n}=1$  and $\d_{n,m}= 0$  if $n\neq m$.

\section{The sigma-function }  

Here we give the precise definition of the sigma-function $\sigma(\lambda)$ and prove the main  identity \e{eq:MAIDs}.  

\medskip

{\bf 2.1.}
We work on test functions $f  \in C_{0}^\infty ({\Bbb R}_{+})$ and   require that $h $ belong to the dual space $ C_{0}^\infty ({\Bbb R}_{+})'$. Let
the set ${\cal Y}$ consist of entire functions  $\varphi(\lambda)$ satisfying,   for all $\lambda \in \Bbb C$,   bounds 
  \begin{equation}
  | \varphi (\lambda)| \leq C_{n}  (1+| \lambda |)^{-n} e^{r_{\pm} |\Re \lambda |}, \q \pm\Re\lambda\geq 0,
 \label{eq:YY}\end{equation}
  for all $n$ and some $r_{+}=r_{+}(\varphi)<0$; the number $r_{-}=r_{-}(\varphi)$  may be arbitrary. The space $\cal Y$ is of course invariant with respect to the complex conjugation  $\varphi(\lambda)\mapsto \varphi^* (\lambda)= \ov{\varphi (\bar{\lambda})}$. By definition, $\varphi_{k} (\lambda)\to 0$ as $k\to \infty$  in $\cal Y$ if all functions $\varphi_{k} (\lambda)$ satisfy bounds \e{eq:YY} with the same constants $r_{\pm}$, $C_{n}$ and $\varphi_{k} (\lambda)\to 0$ as $k\to \infty$  uniformly on all compact subsets of $\Bbb C$. 
  
Let the Laplace  transform  ${\sf L}$ be defined by formula \e{eq:LAPj}.
By one of the versions of the Paley-Wiener theorem,
  ${\sf L}: C_{0}^\infty ({\Bbb R}_{+}) \to {\cal Y}$ is the one-to-one continuous mapping of $ C_{0}^\infty ({\Bbb R}_{+}) $ onto $ {\cal Y}$ and the inverse mapping 
   ${\sf L}^{-1}: {\cal Y} \to  C_{0}^\infty ({\Bbb R}_{+}) $ is also continuous. 
    Passing to the dual spaces, we see that the mapping
      \begin{equation}
      {\sf L}^*: {\cal Y}' \to C_{0}^\infty    ({\Bbb R}_{+})'
 \label{eq:YY1}\end{equation}
 is also an isomorphism. We emphasize that we write ${\sf L} ^*$ here because this operator acts in the spaces of distributions.

   Let us construct the sigma-function. 
   
    \begin{definition}\label{sigmay}
      Assume that 
   \begin{equation}
h \in C_{0}^\infty ({\Bbb R}_{+})' .
 \label{eq:hhg}\end{equation}
   Then the distribution $\sigma \in {\cal Y}'  $ defined by the formula
 \begin{equation}
\sigma=({\sf L}^*)^{-1}h  
 \label{eq:LAPL1g}\end{equation}
is called the sigma-function of the kernel $h$ or of the corresponding Hankel operator $H$.
 \end{definition}

 Since mapping \e{eq:YY1} is an   isomorphism,
  the kernel  $h(t)$ can be recovered
from  its sigma-function $\sigma(\lambda)$ by the formula $h={\sf L}^*\sigma$ which gives the precise sense to  formal relation \e{eq:conv1}. Thus 
   there is the one-to-one correspondence between kernels $  h\in C_{0}^\infty    ({\Bbb R}_{+})' $ and their sigma-functions $\sigma \in {\cal Y} ' $.

      \medskip

{\bf 2.2.}
Now we are in  a position to check the   identity \e{eq:MAIDs}. 
 The first assertion is quite straightforward. It is a direct consequence of Definition~\ref{sigmay}.

 \begin{proposition}\label{1Fg}
Let assumption \e{eq:hhg} hold, and let $\sigma$ be the corresponding sigma-function. 
    Then      the identity
 \begin{equation}
 { \pmb\la} h , F {\pmb \ra} =   { \pmb\la}  {\sf L}^* \sigma ,   F {\pmb \ra} =   { \la} \sigma ,  {\sf L} F {  \ra}
  \label{eq:MAIDG}\end{equation}
  holds for arbitrary
   $F\in C_{0}^\infty    ({\Bbb R}_{+})$.
 \end{proposition}
 
    Let us introduce the   Laplace convolution
\begin{equation} 
( \bar{f}_{1}\star f_{2})(t)=
\int_{0}^t    \overline{f_{1}(s)} f_{2}(t-s) ds  
 \label{eq:HH1}\end{equation} 
 of   functions $ \bar{f}_{1},  f_{2} \in C_{0}^\infty    ({\Bbb R}_{+})$. Then it formally follows from \e{eq:H1}  that
\begin{equation}
(Hf_{1} ,f_{2})=  {\pmb \la} h, \bar{f}_{1}\star f_{2} {\pmb \ra} 
 \label{eq:HH}\end{equation}
 where we write ${\pmb \la} \cdot, \cdot {\pmb \ra}$ instead of $( \cdot, \cdot )$ because $h$ may be a distribution.   
 Obviously,  for  arbitrary ${f}_{1},  f_{2} \in C_{0}^\infty    ({\Bbb R}_{+})$ we have
   \begin{equation}
{\sf L} (\bar{f}_{1} \star  f_{2} ) ={\sf L} \bar{f}_{1}  {\sf L} f_{2} =({\sf L} f_{1} )^* {\sf L} f_{2}\in  {\cal Y}.
 \label{eq:gg1}\end{equation}

   Now we are in a position to precisely state our main  identity. 

 \begin{theorem}\label{1}
Let   assumption \e{eq:hhg} be satisfied, and let $\sigma \in {\cal Y}'$ be defined by formula
  \e{eq:LAPL1g}.
  Then       the identity
 \begin{equation}
 { \pmb\la} h ,  \bar{f}_{1}\star f_{2} {\pmb \ra} =   { \la} \sigma , (  {\sf L} f_{1} )^*{\sf L} f_{2} { \ra}
  \label{eq:MAID}\end{equation}
  holds for arbitrary
   $f_{1}, f_{2}\in C_{0}^\infty    ({\Bbb R}_{+})$.
 \end{theorem}

 \begin{pf}  
 It suffices to apply identity \e{eq:MAIDG} to $F= \bar{f}_{1}\star f_{2}$ and to use relation 
\e{eq:gg1}. 
  \end{pf}

 The identity  \e{eq:MAID} attributes a precise meaning to   \e{eq:MAIDp} or \e{eq:MAIDs}. 
 
 \medskip
 
  {\bf 2.3.}
 Suppose  now that $h(t)=\overline{h(t)}$ for all $t>0$, or to be more precise $\ov{{\pmb \la} h, F{\pmb \ra}} = {\pmb \la} h, \overline{F}{\pmb \ra}$ for all $F\in C_{0}^\infty    ({\Bbb R}_{+})$. Then it follows from  \e{eq:MAIDG} that  the sigma-function is also real, that is,  $\ov{{  \la} \sigma, w{  \ra}} ={  \la} \sigma, w^*{  \ra}$ for all $w\in  {\cal Y} $.

 Below we use the following natural definition.  
 
  \begin{definition}\label{hss}
     Let ${\sf h}[\varphi , \varphi ]$ be   a real quadratic form defined on  a linear set ${\sf D} $.  We denote by $N_{\pm}({\sf h}) = N_{\pm}({\sf h}; {\sf D})$ the maximal dimension of linear sets ${\cal M}_{\pm}\subset {\sf D}$   such that $\pm {\sf h} [\varphi,\varphi] > 0$   for all $\varphi\in {\cal  M}_{\pm}$, $\varphi\neq 0$.
          \end{definition}
          
     Definition~\ref{hss} means that there exists a linear set ${\cal M}_{\pm} \subset {\sf D} $, $\dim {\cal M}_{\pm}= N_{\pm}({\sf h}; {\sf D}) $,
       such that $\pm {\sf h}[\varphi , \varphi ]  > 0$   for all $\varphi\in {\cal M}_{\pm}$, $\varphi \neq 0$, and for every  linear set ${\cal M}_{\pm}' \subset {\sf D}$ with  $\dim {\cal M}_{\pm}'> N_{\pm}({\sf h}; {\sf D})$ there exists $\varphi\in {\cal M}_{\pm}'$, $\varphi \neq 0$,  such that $\pm {\sf h}[\varphi , \varphi ]  \leq 0$.

  Of course, if the set ${\sf D} $ is dense in a Hilbert space $\cal H$ and ${\sf h}[\varphi , \varphi ]$ is semibounded and closed on ${\sf D} $, then for the self-adjoint operator ${\sf  H}$ corresponding to ${\sf h}$, we have  $ N_{\pm}({\sf  H})=N_{\pm}({\sf h}; {\sf D}) $. In particular, this is true for bounded operators ${\sf  H}$.
  
 We apply Definition~\ref{hss} to the forms $h[f,f]={\pmb\la} h, \bar{f} \star f {\pmb\ra}$ on $f \in C_{0}^\infty    ({\Bbb R}_{+})$ and $\sigma[w,w]=\la \sigma, w^* w \ra$  on $w\in   {\cal Y}$.

Since ${\sf L} : C_{0}^\infty    ({\Bbb R}_{+}) \to    {\cal Y}$ is an isomorphism,
   the following assertion is a direct consequence of Theorem~\ref{1}.
 
  \begin{theorem}\label{HBx}
  Let $h \in   C_{0}^\infty    ({\Bbb R}_{+})'$. Then $\sigma = ({\sf L}^*)^{-1}h \in  {\cal Y}' $ and
  \begin{equation}
  N_{\pm}(h; C_{0}^\infty    ({\Bbb R}_{+}))= N_{\pm}(\sigma; {\cal Y} ).
  \label{eq:YS}\end{equation}
   In particular,     the form $\pm {\pmb\la}h,\bar{f} \star f {\pmb\ra} \geq 0$ for all
 $f\in C_{0}^\infty    ({\Bbb R}_{+})$  if and only if the form  $\pm {\la}\sigma , w^* w  {\ra} \geq 0$ for all $w \in {\cal Y} $. 
     \end{theorem}
     
Thus a Hankel operator $H$ is positive (or negative) if and only if its sigma-function $\sigma (\lambda)$ is positive (or negative).

\medskip
 
 {\bf 2.4.}
 In our examples $h(t)$ is a continuous function of $t>0$. However its behavior as $t\to \infty$ and $t\to 0$ may be arbitrary.  

 \begin{example}\label{exp}
  Let $h (t)= e^{t^2}$. Then  
 representation \e{eq:conv1}  is satisfied with the  function $\sigma (\lambda)= 2^{-1}\pi^{-1/2}  e^{-\lambda^2/4}$. Thus ${\pmb \la} h, \bar{f} \star f{\pmb \ra}\geq 0$ for all $f \in   C_{0}^\infty    ({\Bbb R}_{+})$. This result can be compared with the fact that the compact Hankel operator $H$ with kernel $h (t)= e^{-t^2}$ has infinite number of both positive and negative eigenvalues (see Proposition~B.1 in \cite{Y}).    
  \end{example}
  
  In the case considered,  $\supp \sigma={\Bbb R}$ which is by no means true in the general case. For example, if $h(t)=e^{-\alpha t}$ for some $\alpha\in{\Bbb C}$, then 
 $ \la \sigma, w\ra= \ov{w(\alpha)}$. Let us mention a particularly simple special case when the relation $h={\sf L}^* \sigma$ can be understood in the classical sense.
 
  \begin{proposition}\label{HBs}
  Let $h \in   C_{0}^\infty    ({\Bbb R}_{+})'$. With respect to the corresponding sigma-function, assume that
   \begin{equation}
\supp \sigma\subset [0,\infty),
  \label{eq:su}\end{equation}
 $\sigma\in L^1 (0,R)$ for all $R<\infty$ and $\sigma(\lambda)= O (e^{\varepsilon \lambda})$ as $\lambda\to\infty$ for all $\varepsilon>0$. Then
 \[
  h (t)  =    \int_0^\infty     e^{-t \lambda} \sigma(\lambda)  d \lambda.
  \]
   In particular, the function  $h(t)$ is analytic in the right-half plane.
     \end{proposition}

\section{Bounded Hankel operators}  

The main identity \e{eq:MAID}  directly yields a criterion for a Hankel operator to be bounded which provides a new approach to the classical Nehari theorem.

\medskip

 {\bf 3.1.}
 Let $  {\Bbb H}_{r}^p$, $p\geq 1$, be the Hardy space of functions analytic in the right half-plane. 
 Obviously, $w \in  {\Bbb H}_{r}^2$ if and only if its complex conjugate $w^* \in  {\Bbb H}_{r}^2$.
         By the Paley-Wiener theorem, the operators    
           \[
           (2\pi)^{-1/2} {\sf L}: L^2 ({\Bbb R}_{+}) \to {\Bbb H}_{r}^2 \q \mbox{and hence}\q 
  (2\pi)^{-1/2} {\sf L}^*: ({\Bbb H}_{r}^2)' \to L^2 ({\Bbb R}_{+}) 
\]
are  unitary.
         Moreover, $\| {\sf L} f\|_{L^2 ({\Bbb R}_{+})} \leq \sqrt{\pi} \|   f\|_{L^2 ({\Bbb R}_{+})}$ (see, e.g., formulas \e{eq:LAPL1} and \e{eq:sig1x} below). Putting $w={\sf L}f$, we see that
             \[
\sqrt{2} \| w \|_{L^2 ({\Bbb R}_{+})} \leq \| w \|_{{\Bbb H}_{r}^2}.
\]

Let us state a criterion for  boundedness of a Hankel operator $H$ in terms of its sigma-function.
We proceed from the definition  of   $H$ by  its quadratic form \e{eq:HH}  where  $
f_{j}\in C_{0}^\infty ({\Bbb R}_{+})$, $j=1,2$, and  $h\in C_{0}^\infty ({\Bbb R}_{+})'$. Recall that, by Definition~\ref{sigmay}, in this case   its sigma-function $\sigma\in {\cal Y}'$. Of course $  {\cal Y}\subset {\Bbb H}_{r}^1$ so that the dual space
$  ({\Bbb H}_{r}^1)' \subset {\cal Y}'$.
       
         \begin{theorem}\label{NehC}
        A Hankel operator $H$ is bounded  in the space $L^2 ({\Bbb R}_{+})$   if and only if its sigma-function $\sigma\in ({\Bbb H}_{r}^1)'$.  
   \end{theorem}
   
   \begin{pf}
  According to the identity \e{eq:MAID} 
   $H$ is bounded   if and only if 
     \[
|\la \sigma,  ({\sf L} f_{1})^*{\sf L} f_{2}\ra |\leq C \|   f_{1}\|_{L^2 ({\Bbb R}_{+})} \|   f_2\|_{L^2 ({\Bbb R}_{+})}
\]
for all  $f_{j}\in C_{0}^\infty ({\Bbb R}_{+})$ or, equivalently, all $f_{j}\in L^2 ({\Bbb R}_{+})$. Putting $w_{j}= {\sf L} f_j$, we
 rewrite this estimate   as 
 \begin{equation}
|\la \sigma,  w_{1}^* w_{2}\ra |\leq C \| w_{1}\|_{{\Bbb H}_{r}^2} \| w_2\|_{{\Bbb H}_{r}^2},\q \forall w_{j}\in {\Bbb H}_{r}^2.
\label{eq:NES1}\end{equation}
It is obviously satisfied if $\sigma\in ({\Bbb H}_{r}^1)'$. 

 Conversely, in view of the inner-outer factorization (see, e.g., \cite{NK}), every $g\in{\Bbb H}_{r}^1$ admits the representation $g=w_{1}^*w_{2}$ where 
$w_{1}, w_{2}\in{\Bbb H}_{r}^2 $ and
\[
\|  g    \|_{{\Bbb H}_{r}^1} =\|  w_1    \|_{{\Bbb H}_{r}^2} \| w_2   \|_{{\Bbb H}_{r}^2}.
\]
Therefore according to \e{eq:NES1} we have
 \[
| {\la } \sigma, g {\ra } | \leq C \|  g    \|_{{\Bbb H}_{r}^1}  , \q \forall g\in {\Bbb H}_{r}^1, 
 \]
whence $\sigma\in ({\Bbb H}_{r}^1)'$. 
      \end{pf}

       \medskip
  
   {\bf 3.2.}
   Theorem~\ref{NehC}  can equivalently be reformulated in terms of symbols of Hankel operators. This requires the Fefferman duality result (see the original paper \cite{Fef} or Theorem~4.4 in Chapter~VI of the book \cite{Garnett}). We denote by ${\sf B}_{r}$ the class of analytic in the right half-plane functions which have a bounded mean oscillation on the imaginary axis. We omit standard explanations of the precise meaning of the 
   integral in the right-hand side of \e{eq:Fef}.
   
    \begin{theorem}[Fefferman]\label{NeFF}
    A functional    $\sigma\in ({\Bbb H}_{r}^1)'$ if and only if there exists a function $\omega \in{\sf B}_{r}$ such that
     \begin{equation}
 {\la } \sigma, g {\ra } = -i \int_{-i\infty}^ {i\infty }\omega(\mu) \ov{g(\mu)} d\mu
 \label{eq:Fef}\end{equation}
 for all    $g\in {\Bbb H}_{r}^1$.
   \end{theorem}
   
   Now it is easy to deduce the  classical Nehari-Fefferman result  from  Theorem~\ref{NehC}. We recall that the symbol $\omega$ of a Hankel operator $H$ is defined by formula \e{eq:hs2b} so that
   \begin{equation}
 {\pmb\la } h, \bar{f}_{1}  \star f_{2}   {\pmb\ra }    = -i  \int_{-i\infty}^{i\infty} \omega(\mu) 
     \ti{f}_{1} (-\mu) \ov{\ti{f}_{2} ( \mu)}d\mu, \q \ti{f}_j= {\sf L}f_{j}.
    \label{eq:FEH}\end{equation} 
   
     \begin{theorem}\label{NehCF}
        A Hankel operator $H$ is bounded  in the space $L^2 ({\Bbb R}_{+})$   if and only if its symbol $\omega \in {\sf B}_{r}$.    
   \end{theorem}
   
   \begin{pf}
   If $H$ is bounded, then, by Theorem~\ref{NehC}, its sigma-function $\sigma\in ({\Bbb H}_{r}^1)'$. By Theorem~\ref{NeFF},  relation  \e{eq:Fef} is satisfied with some $\omega \in{\sf B}_{r}$. Now 
using  the main identity  \e{eq:MAID} and
  applying relation  \e{eq:Fef} to   the  function $g=({\sf L}f_1)^* {\sf L}f_2$, we get \e{eq:FEH}.
  
  Conversely, let $\omega \in{\sf B}_{r}$. Then according to Theorem~\ref{NeFF} it follows from \e{eq:FEH} that
  \[
|   {\pmb\la }  h, \bar{f}_{1}  \star f_{2}  {\pmb\ra }  | \leq C \|   \ti{f}_{1}^*      \ti{f}_2   \|_{{\Bbb H}_{r}^1}
\leq   C  \|   \ti{f}_1    \|_{{\Bbb H}_{r}^2} \|  \ti{f}_2   \|_{{\Bbb H}_{r}^2}= 2\pi
C \|  f_1    \|_{L^2 ({\Bbb R}_{+})} \| f_2   \|_{L^2 ({\Bbb R}_{+})}.
\]
Thus $H$ is bounded.
   \end{pf}

 We emphasize that in contrast to the original proof of the Nehari theorem (see his paper \cite{Nehari} or the book 
 \cite{Pe}), the proof of Theorem~\ref{NehC}  does not require either the Hahn-Banach or M.~Riesz  theorems. Only the inner-outer factorization has been used.

   \medskip
  
   {\bf 3.3.}
   Let us illustrate the link between $\sigma$ and $\omega$ at the example of the Carleman operator with kernel $h(t)=t^{-1}$. According to \e{eq:conv1} we have $\sigma(\lambda)=1$ for $\lambda\in {\Bbb R}_{+}$ and $\sigma(\lambda)=0$ for $\lambda \not \in {\Bbb R}_{+}$. 
   
   Let us show that $\omega (\mu)=  -\ln\mu$.   According to formula \e{eq:hs2a} we only have to check that
  \begin{equation}
 2\pi i  \int_{0}^{\infty}   \ov{g(\lambda)}d\lambda  = - \int_{-i\infty}^{i\infty} \ln\mu \,\ov{g(\mu)}d\mu 
    \label{eq:hs2C}\end{equation}  
    for $g\in {\Bbb H}_{r}^1$. It suffices to consider the functions $g(\mu)=(\mu+a)^{-n}$ for $n\in{\Bbb Z}_{+}$, $n\geq 2$, $a>0$. The right-hand side of \e{eq:hs2C} equals
    \[
    - \int_{-i\infty}^{i\infty} \ln\mu \, (-\mu+ a)^{-n} d\mu  = 2\pi i (-1)^n \Res_{\mu=a} \big(\ln\mu \,(\mu-a)^{-n}\big)=
    2\pi i (n-1)^{-1} a^{-n+1}
    \]
   which obviously coincides with the left-hand side of   \e{eq:hs2C}.

  Alternatively, for the calculation of $w(\mu)$, we can proceed from
    Theorem~8.8 of Chapter~1 in the book \cite{Pe}. To that end, we first  have  to extend the distribution $h(t)=t^{-1}$ from $C_{0}^\infty ({\Bbb R}_{+})$   onto the Schwartz space ${\cal S} ({\Bbb R}_{+})$. This is done by the formula
   \[
{\pmb\la } h, \varphi {\pmb\ra }= \int_{0} ^1 \frac{\bar{\varphi}(t)-\bar{\varphi}(0) } {t}dt
 +  \int_1^\infty \frac{\bar{\varphi}(t)  } {t}dt.
 \]
 Therefore according to formula \e{eq:hst} we have
\[
\omega(\mu)= \int_{0} ^1 \frac{e^{-\mu t}-1 } {t}dt
 +  \int_1^\infty \frac{e^{-\mu t} } {t}dt= -\ln\mu+ \Gamma' (1).
 \]
 The constant term here can be of course neglected.

\section{A factorization of the Laplace transform}


In this section we collect various results which rely on a factorization of the Laplace transform.

\medskip
 
 {\bf 4.1.}
   For a factorization of the Laplace transform $\sf L$, it is natural to consider more general integral operators 
   \begin{equation}
(A f)(t)= \int_{0}^\infty a(ts) f(s) ds
 \label{eq:ME}\end{equation}
 with kernels $a$ depending on the product of the variables only. Such operators   can be standardly diagonalized (see, e.g., \cite{Y3}) by the Mellin transform $M $. Let   a unitary operator $U: L^2 ({\Bbb R}_{+})\to L^2 ({\Bbb R} )$ be defined by the formula 
  \begin{equation}
(U f)(x) =e^{x/2} f(e^x).
 \label{eq:HHU}\end{equation}
 Then $M =\Phi U$.
 
 We suppose that the function $a(t) t^{-1/2  } $ belongs to $L^1 ({\Bbb R}_{+})$ (in this case the operator $A$ is bounded in the space $L^2 ({\Bbb R}_{+})$). 
Making in \e{eq:ME} the change of variables $t=e^x$, $s=e^y$, we see  that 
   \[
(UA f)(x)= \int_{-\infty}^\infty  (Ua) (x+y) (Uf)(y) dy, \q f\in L^2 ({\Bbb R}_{+}).
\]
 Passing here to the Fourier transforms, we find that
  \begin{equation}
( M A f)(\xi)=  \sqrt{2\pi} (Ma ) (\xi) (Mf)(-\xi), \q M=\Phi A .  
 \label{eq:ME3}\end{equation}
Let  ${\cal J}$, $({\cal J}u)(\xi)= u(-\xi)$, be the reflection operator and 
   \begin{equation}
   {\sf a} (\xi)= \sqrt{2\pi} (M a)(-\xi)  =\int_{0}^\infty a(t) t^{-1/2 + i\xi} dt.
 \label{eq:ME4}\end{equation}
It follows from \e{eq:ME3} that
   \begin{equation}
 A f = M^{-1} {\cal J}{\sf a} M f. 
 \label{eq:ME5}\end{equation}
 Thus we obtain the following assertion.
 
   \begin{theorem}\label{ME}
   Let the operator $A$ be defined by formula \e{eq:ME}  where
    the function $a(t) t^{  -1 /2  } $ belongs to $L^1 ({\Bbb R}_{+})$, and let $ {\sf a} (\xi)$ be defined by formula \e{eq:ME4}. 
    Then for all  $ f\in L^2 ({\Bbb R}_{+})$ representation \e{eq:ME5} holds.
    \end{theorem}

  Let us apply this result to the case $a(t)=e^{-t}$ when $A={\sf L}$ and
    \[
   {\sf a}  (\xi) =\int_{0}^\infty e^{-t} t^{ -1/2 + i\xi} dt=\Gamma (1/2 +i\xi)
\]
 is the gamma function.  Recall that the gamma function $\Gamma (z)\neq 0$ for all $z\in {\Bbb C}$. According to the Stirling formula, the function
\begin{equation}
   \Gamma(\gamma +i \xi) = e^{\pi i (2\gamma-1)/4} (2\pi / e)^{1/2}  \xi^{\gamma-1/2} e^{i\xi(\ln \xi -1)} e^{-\pi \xi /2}\big(1+O(\xi^{-1})\big)  
\label{eq:M11}\end{equation}
  tends exponentially to zero as $\xi\to +\infty$ if $\gamma>0$ is fixed.  Since $\Gamma (\gamma- i\xi) =\ov{\Gamma (\gamma+i\xi)}$, the same is true as $\xi\to -\infty$. We put
   \[
  ({\pmb \Gamma}_{\gamma}   u)(\xi)=\Gamma(\gamma +i\xi) u(\xi). 
 \]
 
 Theorem~\ref{ME} implies the following statement.

  \begin{corollary}\label{LAPL} 
  For all $f\in L^2 ({\Bbb R}_{+})$,  the identity
    \begin{equation}
({\sf L}f )(\lambda)=  (M^{-1} {\cal J} {\pmb \Gamma}_{1/2}  M  f )(\lambda), \q \lambda>0,
 \label{eq:LAPL1}\end{equation}
 holds.
  \end{corollary}   

   This formula   can   be used for the inversion of the Laplace transform:
   \[
{\sf L}^{-1}=  M^{-1}   {\pmb \Gamma}_{1/2}^{-1} {\cal J}  M  .
 \]
 Observe   that according to \e{eq:M11} the function $\Gamma(1/2 +i\xi)^{-1} $   exponentially grows as $|\xi|\to\infty$.  This is why,  even for very nice kernels $h(t)$ (for example, for   $h(t)= t^k e^{- \alpha t}$, $\Re\alpha >0$, $k=0,1,\ldots$), the corresponding sigma-function $\sigma (\lambda)$ defined by formula \e{eq:ss} may be a highly singular distribution.

 \medskip
 
 {\bf 4.2.} 
Factorization \e{eq:LAPL1} allows us to reformulate the main identity \e{eq:MAID} in a somewhat  different form. We suppose for simplicity that conditions \e{eq:su}  and
 \begin{equation}
\sigma \in L^\infty ({\Bbb R}_{+})
 \label{eq:sig1}\end{equation}
are satisfied. Then  the operators $\Sigma$ and hence $H$ are bounded  in the space $L^2 ({\Bbb R}_{+})$.  Let the function $s(x)$ be  defined  for $x\in {\Bbb R}$ by formula
 \e{eq:sig2}, and let $S$ be the operator of multiplication by   $s(x)$ in the space $L^2 ({\Bbb R} )$.
    Since $M=\Phi U$, we have
  \begin{equation}
  {\cal J}  M \Sigma M^{-1} {\cal J} =\Phi S \Phi^*.
 \label{eq:sig3}\end{equation}

 Let us further observe that
 \begin{equation}
 |\Gamma (1/2+ i\xi)|= \frac{\sqrt{\pi} }{\sqrt{\cosh(\pi\xi)} }=: v(\xi)
 \label{eq:sig1x}\end{equation}
 and   denote by $V$ the operator of multiplication by this function in the space $L^2 ({\Bbb R} )$. Set also
 \[
 ({\bf M} f) (\xi)=e^{i\arg \Gamma (1/2+ i\xi)}  (M f) (\xi).
\]

Putting together the identities \e{eq:MAIDp}, \e{eq:LAPL1}  and \e{eq:sig3}, we obtain the following result.

   \begin{theorem}\label{SIG} 
 Under assumptions  \e{eq:su}  and \e{eq:sig1} define $s(x)$ by formula \e{eq:sig2} and set
   \[
  A = V \Phi S \Phi^* V.
 \]
 Then
  \begin{equation}
H= {\bf M}^* A {\bf M}.
 \label{eq:sig6}\end{equation}
   \end{theorem}  
   
Note that according to formula \e{eq:conv1} under the assumptions of this theorem the kernel $h(t)$ of $H$ satisfies the bounds
   \[
 |  h^{(n)}(t) | \leq C_{n} t^{-1-n}, \q \forall n\in{\Bbb Z}_{+}.
   \]
   
   Obviously, $A$ is a pseudo-differential operator in the space $L^2 ({\Bbb R})$ with  
   amplitude  \e{eq:sam} which  factorizes into a product of functions of   one variable only. Assumption \e{eq:sig1} is by no means necessary. For example, if $h(t)= P(\ln t) t^{-1}$ where  $P$ is   a polynomial, then $s(x)$ is also a polynomial (see \cite{Yf1}, for details). In this case relation \e{eq:sig6}  holds with   a differential operator $A$.
      Thus   Theorem~\ref{SIG} shows that, under very general assumptions,    Hankel operators are unitarily equivalent to pseudo-differential operators  of   a very special structure.

 \medskip
 
 {\bf 4.3.} 
 Let us come back to Theorem~\ref{HBx}.
 It is usually not   convenient to work with analytic test functions. Fortunately under mild additional assumptions on the sigma-function, the set $\cal Y$  of test functions $w$ in \e{eq:YS} can be replaced by the set $C_{0}^\infty    ({\Bbb R}_{+})$.
  Below we sometimes do not distinguish functions $g\in \cal Y$ and their restrictions on ${\Bbb R}_{+}$.

   Let us  introduce the  space ${\cal S}_\gamma$, $\gamma \in {\Bbb R}$,   of functions $g\in C^\infty ({\Bbb R}_{+})$ satisfying estimates 
    \begin{equation}
| g^{(k)} (\lambda)| \leq C_{\varkappa,k} \lambda^{\gamma -1/2 -k} (1+| \ln \lambda |)^{-\varkappa}
 \label{eq:SZ}\end{equation} 
 for all  $k=0,1,2,\ldots$ and all $\varkappa\in {\Bbb R}$. The case $\gamma=0$ is the most important for us. It is easy to see that $g\in {\cal S}_0$ if and only if the
 function $ Ug $ belongs to the Schwartz space ${\cal S}$.

  We need the following analytical result.
 
  \begin{lemma}\label{SZ}
  The set $\cal Y$ is dense in ${\cal S}_0$.
     \end{lemma}
  
   \begin{pf}
   The result formulated is equivalent to the fact that the set of elements $U {\sf L} f$ where $f \in C_{0}^\infty ({\Bbb R}_{+} )$ is dense in the   space ${\cal S}$.
      In view of the identity \e{eq:LAPL1}, it is equivalent to the following assertions:  the set of elements $\Phi^{-1} J {\pmb \Gamma}_{1/2}  M f$ where $f \in C_{0}^\infty ({\Bbb R}_{+} )$ is dense in the   space ${\cal S}$ or the set of elements $  {\pmb \Gamma}_{1/2}  \Phi\psi $ where $\psi= U f \in C_{0}^\infty ({\Bbb R} )$ is dense in the   space ${\cal S}$.

    Thus, for an arbitrary $u\in {\cal S}$, we have to construct a sequence $\psi_{k} \in C_{0}^\infty ({\Bbb R} )$ such that
    \begin{equation}
{\pmb \Gamma}_{1/2}   \Phi \psi_{k}  \to u
 \label{eq:SZ1}\end{equation} 
 in ${\cal S}$ as $k\to \infty$.
   Let $\theta \in C_{0}^\infty ({\Bbb R} )$ and $\theta(\xi)=1$ in a neighborhood of the point $\xi=0$. Put $\theta_n(\xi)= \theta(  \xi/n)$. Then
    \begin{equation}
 \theta_{n} u \to u 
 \label{eq:SZ2}\end{equation} 
 in ${\cal S}$ as $n\to \infty$. 
 Next, we put  
     \[
v_n (\xi)= \Gamma(1/2 +i\xi)^{-1} \theta_n (\xi) u(\xi).
 \]
 Obviously,  $v_n \in C_{0}^\infty ({\Bbb R} ) \subset {\cal S}$ and hence, for every $n$, there exists a sequence $\psi_{n,m} \in C_{0}^\infty ({\Bbb R} )$ such that   $\Phi \psi_{n,m}  \to v_n$  in ${\cal S}$. It follows that
  \begin{equation}
{\pmb \Gamma}_{1/2}     \Phi \psi_{n,m}  \to {\pmb \Gamma}_{1/2}   v_n =\theta_n u
 \label{eq:SZ4}\end{equation} 
 in ${\cal S}$ as $m\to \infty$.  Putting together relations \e{eq:SZ2} and \e{eq:SZ4}, we can choose  a subsequence $\psi_{k}  $ of $\psi_{n,m} $ such that relation \e{eq:SZ1} is true.
       \end{pf}

  Lemma~\ref{SZ} allows us to prove the following assertion.

  \begin{lemma}\label{HBxe}
  Under assumption \e{eq:su} 
suppose   that $\sigma   \in  {\cal S}'_{  0}$. 
  Then     
   \begin{equation}
    N_{\pm}(\sigma; {\cal Y} )=N_{\pm}(\sigma; {\cal S}_{0} )=N_{\pm}(\sigma; C_{0}^\infty ({\Bbb R}_{+}) ).
\label{eq:cb}\end{equation} 
  \end{lemma}
  
  \begin{pf}
  Let us check the first equality \e{eq:cb}.
  The inequality 
$ N_{\pm}(\sigma; {\cal Y} )\leq N_{\pm}(\sigma; {\cal S}_{0} ) $  is obvious because $\cal Y\subset {\cal S}_{0}$.     
  
      Let us prove the opposite inequality. Consider for definiteness the sign $``+"$. Let ${\cal L} \subset {\cal S}_{0}$, and let $  \sigma[w,w]  > 0$   for all $w\in {\cal L}$, $w \neq 0$.   
    Suppose first that $N : = \dim {\cal L} < \infty$ and choose elements $w_1,\ldots ,w_{{N}} \in {\cal L}$ such that $\sigma[w_j, w_k ]=\delta_{j,k}$ for all $j,k=1,\ldots, {N}$. Using Lemma~\ref{SZ} we can construct elements $w_j^{(\epsilon)}  \in  {\cal Y}$ 
   such that $w_j^{(\epsilon)} \to w_j$ and hence $w_j^{(\epsilon)} \bar{w}_k^{(\epsilon)} \to w_j \bar{w}_k$ in ${\cal S}_{0}$ as $\epsilon\to 0$ for all $j,k=1,\ldots,  {N}$.  Since $\sigma\in {\cal S}'_0$, we see that $\sigma[ w_j^{(\epsilon)} , w_k^{(\epsilon)} ]\to \delta_{j,k}$ as $\epsilon\to 0$. For an arbitrary $ \gamma >0$, we can choose $\epsilon$ such that 
  $| \sigma[ w_j^{(\epsilon)} , w_k^{(\epsilon)} ] - \delta_{j,k}| \leq \gamma  $. Then for arbitrary $\lambda_1, \ldots, \lambda_{{N}}\in {\Bbb C}$, we have
   \begin{multline}
\sigma [\sum_{j=1}^{{N}}\lambda_j w_j^{(\epsilon)} , \sum_{j=1}^{{N}}\lambda_j w_j^{(\epsilon)} ]=
\sum_{j=1}^{{N}} |\lambda_j |^2 \sigma [w_j^{(\epsilon)}, w_j^{(\epsilon)}]+ 2\Re \sum_{j,k=1; j\neq k}^{{N}} 
\lambda_j \bar{\lambda_k} \sigma [w_j^{(\epsilon)}, w_k^{(\epsilon)}]
 \\
 \geq (1-\gamma )\sum_{j=1}^{{N}} |\lambda_j |^2  -2 \gamma  \sum_{j,k=1; j\neq k}^{{N}} 
\lambda_j \bar{\lambda_k}
 \geq \big(1-(2 {N}-1)\gamma \big)\sum_{j=1}^{{N}} |\lambda_j |^2 . 
  \label{eq:abcs}\end{multline}
  Thus elements $w_1^{(\epsilon)}, \ldots, w^{(\epsilon)}_{{N}}$ are linearly independent if
  $(2 {N}-1)\gamma < 1$.  The same inequality \e{eq:abcs} shows that $\sigma[w,w]>0$ on all vectors $w\neq 0$ belonging to the space $ {\cal L}^{(\epsilon)} $ spanned by $ w_1^{(\epsilon)}, \ldots,  w_ {{N}}^{(\epsilon)}  $.

  If $ {N} =\infty$, then the  construction above works on every finite dimensional subspace of ${\cal L} $ where  $\sigma[w,w]>0$. This yields the  space ${\cal L}^{(\epsilon)}   \subset {\cal Z}$ of an arbitrary large dimension where $\sigma[w,w]>0$.
  
  The second equality \e{eq:cb} can be proven quite similarly because $C_{0}^\infty ({\Bbb R}_{+}) \subset {\cal S}_{0}$, and it is dense in ${\cal S}_{0}$.
   \end{pf}

Putting together Theorem~\ref{HBx} and Lemma~\ref{HBxe} we arrive at the following result.

    \begin{theorem}\label{HByz}
 Let $h \in  C_{0}^\infty ({\Bbb R}_{+})'$.  Suppose   that condition \e{eq:su} is satisfied and that $\sigma   \in  {\cal S}'_{ 0}$. Then 
    \begin{equation}
   N_{\pm}(h; C_{0}^\infty ({\Bbb R}_{+}) ) =N_{\pm}(\sigma; C_{0}^\infty ({\Bbb R}_{+})).
\label{eq:abc}\end{equation}
   \end{theorem}

     The following consequence of Theorem~\ref{HByz} is very convenient for applications to concrete  Hankel operators.

      \begin{theorem}\label{HBy}
  Let the assumptions of Theorem~\ref{HByz} hold. 
  
  $1^0$  If $\pm\sigma\geq 0$, then $\pm{\pmb\la}h, \bar{f}\star f {\pmb\ra} \geq 0$ for all $f\in C_{0}^\infty ({\Bbb R}_{+})$.
  
  $2^0$ If  $\sigma\in L^1_{\rm loc}(\Delta)$ for some interval $\Delta\subset{\Bbb R}$ and
  $\pm\sigma (\lambda)\geq \sigma_{0}>0$  for almost all 
  $\lambda\in \Delta $, then $   N_{\pm}(h; C_{0}^\infty ({\Bbb R}_{+}) ) =\infty$.
   \end{theorem}
   
   \begin{pf}
   The first assertion is a direct consequence of relation \e{eq:abc}.  
   
   For the proof of the second     assertion, choose some
 number $N$ and a function $\varphi\in C_{0}^\infty ({\Bbb R}_{+})$ such that $\varphi(\lambda)=1$ for $\lambda\in [-\d,\d]$ and $\varphi(\lambda)=0$ for $\lambda\not\in [-2\d,2\d]$ where $\d=\d_{N}$ is a sufficiently small number. Let points $\alpha_{j}\in \Delta$, $j=1,\ldots, N$, be such that
   $\alpha_{j+1}-\alpha_{j} =\alpha_{j }-\alpha_{j-1}$ for $j=2,\ldots, N -1$. Set $\Delta_{j}= (\alpha_{j}-\d, \alpha_{j} +\d)$, $\wt{\Delta}_{j}= (\alpha_{j}-2\d, \alpha_{j} +2\d)$. For a sufficiently small $\d$, we may suppose that $\wt{\Delta}_{j}\subset \Delta$ for all $j=1,\ldots, N$ and that $\wt{\Delta}_{j}\cap \wt{\Delta}_{j+1}=\varnothing$ for   $j=1,\ldots, N -1$. We set $\varphi_{j}(\lambda)=\varphi (\lambda-\alpha_{j})$. Since $\pm\sigma (\lambda)\geq \sigma_{0}>0$   for   $\lambda \in \Delta$, we have   
    \[
\pm \sigma [\varphi_{j} , \varphi_{j}] = \pm\int_{0}^\infty \sigma(\lambda) |\varphi_{j}(\lambda) |^2 d\lambda\geq 2\d \sigma_{0}>0.
\]
 The functions $\varphi_{1},\ldots,\varphi_N$ have disjoint supports and hence  $\pm \sigma [w , w] > 0$ for an arbitrary nontrivial  linear combination $w$ of the functions $\varphi_{j}$. Therefore $N_{\pm}(\sigma; C_{0}^\infty ({\Bbb R}_{+}))\geq N$. Since $N$ is arbitrary,  it remains   to use  relation \e{eq:abc}.
        \end{pf}

\medskip

{\bf 4.4.}  
Here we define a scale of spaces of analytic functions where the Laplace transform acts as an isomorphism.  This extends the one-to-one correspondence of Section~2 between kernels of Hankel operators and their sigma-functions to new spaces of distributions.

Let the space $\cal Z= \cal Z ({\Bbb R})$ of test functions be defined as the subset  
 of the Schwartz  space ${\cal S}={\cal S} ({\Bbb R}) $ which consists of functions $\varphi  (x)$ admitting the analytic continuation to   entire functions in the   complex plane $\Bbb C$   and satisfying,   for all $z\in \Bbb C$,   bounds 
  \begin{equation}
  | \varphi (z)| \leq C_{\varkappa}  (1+| z |)^{-\varkappa} e^{r |\Im z |}
 \label{eq:ZZ}\end{equation}
  for some $r=r(\varphi)>0$  and all $\varkappa$. The space $\cal Z$ is of course invariant with respect to the complex conjugation, that is, $\varphi^* (z) = \ov{\varphi (\bar{z})}$ belongs to $\cal Z$  together with $\varphi$. By definition, $\varphi_{k} (z)\to 0$ as $k\to \infty$  in $\cal Z$ if all functions $\varphi_{k} (z)$ satisfy bounds \e{eq:ZZ} with the same constants $r$, $C_{\varkappa}$ and $\varphi_{k} (z)\to 0$ as $k\to \infty$  uniformly on all compact subsets  of $\Bbb C$.
Recall (see, e.g., the book \cite{GUECH})  that the Fourier transform 
  $\Phi $ is a one-to-one mapping of $\cal Z $ onto $C_{0}^\infty ({\Bbb R})$. Moreover, $\Phi$ as well as its inverse   $\Phi^{-1}$ are continuous mappings so  that $\Phi : \cal Z \to C_{0}^\infty ({\Bbb R})$  is an isomorphism.

 Let $U$ be operator \e{eq:HHU}. We define  the set ${\cal Z}_0$ of test functions $f(t)$ by the condition
   \[
f\in  {\cal Z}_0 \Longleftrightarrow Uf\in  {\cal Z}.
\]
 The set ${\cal Z}_0\subset L^2 ({\Bbb R}_{+})$, and it is dense in $L^2 ({\Bbb R}_{+})$ because ${\cal Z}$ is dense in $L^2 ({\Bbb R})$.
 Define  also the set ${\cal Z}_\gamma$ for an arbitrary  $\gamma\in {\Bbb R}$   by the condition $f\in {\cal Z}_\gamma$ if and only if the
 function $t^{-\gamma}f(t)$ belongs to ${\cal Z}_0$.   Thus $f\in {\cal Z}_\gamma$ if and only if the
 function $ F(x)=e^{(1/2-\gamma)x}f(e^x)$ belongs to ${\cal Z}$, that is,
  \[
f (t)= t^{\gamma -1/2} F(\ln t)
\]
 where $F\in  {\cal Z}$. Functions $f(t)$ admit analytic continuation $f(\zeta)$ onto the Riemann surface of the logarithmic function, and   they satisfy the bounds
  \begin{equation}
|f(\zeta)| \leq C_{\varkappa} |\zeta|^{ \gamma-1/2} (1+\big|\ln  |\zeta| \big|)^{-\varkappa}e^{r |\arg \zeta|}
 \label{eq:heta1}\end{equation}
 with some constant $r=r(f)>0$ for all $\varkappa \in {\Bbb R}$. Note that ${\cal Z}_\gamma\subset {\cal S}_\gamma$ where the set ${\cal S}_\gamma$ is defined by conditions  \e{eq:SZ}.
 The sets ${\cal Z}_\gamma$ are invariant with respect to the complex conjugation because ${\cal Z}$  is.
 The topology on  ${\cal Z}_\gamma$  is of course induced by that on ${\cal Z} $. 
 
 Clearly, for all $\gamma,\beta \in {\Bbb R}$, a function  $f\in{\cal Z}_\gamma$ if and only if the function  $t^\beta f(t)$ belongs to ${\cal Z}_{\gamma+\beta}$.  Note that there is no ordering between different spaces ${\cal Z}_\gamma$. If $\gamma_{2}>\gamma_{1}$, then functions $f\in{\cal Z}_{\gamma_{2}}$ are better then those in ${\cal Z}_{\gamma_{1}}$ as $t\to 0$ but worse as $t\to \infty$. Of course neither of the inclusions ${\cal Z}_{\gamma}\subset C_{0}^\infty({\Bbb R}_{+})$ nor $ C_{0}^\infty({\Bbb R}_{+})\subset {\cal Z}_{\gamma} $ (for any $\gamma$) is true.

Since the product of two functions in $\cal Z$ also belongs to this space, the statement below is a direct consequence of the definition of ${\cal Z}_\gamma$.

  \begin{lemma}\label{zz}
If $f\in {\cal Z}_\gamma$ and $g\in {\cal Z}_\beta$ for some $\gamma, \beta\in {\Bbb R}$, then 
$f g \in {\cal Z}_{\gamma+\beta -1/2}$.
  \end{lemma}
  
   Applying now Theorem~\ref{ME} to the kernel  $a(t)=e^{-t}t^{\gamma-1/2}$, we obtain the following generalization of Corollary~\ref{LAPL}.
  
    \begin{lemma}\label{MEan}
    Let $\Omega$ be the operator of multiplication by $t$ or $\lambda$. Then,
    for all $\gamma>0$, the representation
      \[
{\sf L} =\Omega^{1/2 - \gamma}  M^{-1} {\cal J}{\pmb  \Gamma}_{\gamma} M  \Omega^{1/2 - \gamma} 
\]
 holds.
   \end{lemma}

   Observe now that $ U: {\cal Z}_{0} \to {\cal Z}$, $\Phi : {\cal Z} \to C_{0}^\infty({\Bbb R})$ and hence $M: {\cal Z}_{0} \to C_{0}^\infty({\Bbb R})$ are isomorphisms. Moreover, ${\cal J}{\pmb  \Gamma}_{\gamma} : C_{0}^\infty({\Bbb R}) \to C_{0}^\infty({\Bbb R})$ is an isomorphism
because   $\Gamma(\gamma+i\xi) \neq 0$ for all $\xi\in {\Bbb R}$. Therefore Lemma~\ref{MEan} yields

  \begin{corollary}\label{MEanx}
    For all $\gamma>0$,   the mapping
  \begin{equation}
{\sf L} : {\cal Z}_{\gamma-1/2} \to {\cal Z}_{1/2-\gamma}
 \label{eq:LL1}\end{equation}
 is an isomorphism.   
  \end{corollary}

 Note that functions $f(t)$ in ${\cal Z}_{1/2- \gamma}$ are better (worse) than those in ${\cal Z}_{\gamma-1/2}$ at infinity (at zero) if $ \gamma  \leq 1/2$. It is the opposite if $ \gamma \geq 1/2$. For $\gamma=1/2$, the mapping ${\sf L} : {\cal Z}_{0} \to {\cal Z}_{0}$ is an automorphism.

        Let    ${\cal Z}_{\gamma}'$ be the space  dual to ${\cal Z}_{\gamma}$.
        Obviously, $h\in  {\cal Z}_{\gamma}'$ if and only if the function $e^{(\gamma+1/2)x} h(e^x)$ belongs to the space  $   {\cal Z}'$ dual to $   {\cal Z}$. Since ${\cal Z}_{\gamma+\beta}=\Omega^\gamma {\cal Z}_\beta$, we have ${\cal Z}_{\gamma+\beta}'=\Omega^{- \gamma} {\cal Z}_\beta'$ for all $\gamma,\beta\in {\Bbb R}$.  Note that for all $\gamma,\beta \in {\Bbb R}$, a distribution  $h\in{\cal Z}'_\gamma$ if and only if the distribution  $t^\beta h(t)$ belongs to ${\cal Z}'_{\gamma+\beta}$.

 It follows from \e{eq:LL1} that the mapping
  \begin{equation}
{\sf L} ^*: {\cal Z}_{1/2- \gamma}' \to {\cal Z}_{ \gamma-1/2}' ,\q \forall \gamma>0, 
 \label{eq:LL2}\end{equation}
is an isomorphism which according to Definition~\ref{sigmay} yields the following result.

 \begin{proposition}\label{1FD}
For all $\gamma>0$, there is the one-to-one correspondence between kernels $  h\in {\cal Z}_{\gamma-1/2}' $ and their sigma-functions $\sigma \in {\cal Z}_{ 1/2-\gamma}' $, that is, 
 \[
h\in {\cal Z}_{\gamma-1/2}' \Longleftrightarrow\sigma =( {\sf L} ^*)^{-1} h\in {\cal Z}_{ 1/2-\gamma}' 
\]
 \end{proposition}
 
 
  It follows from condition \e{eq:SZ}  that $ h\in {\cal S}_{\gamma-1/2}' \subset {\cal Z}_{\gamma-1/2}' $  if $h\in L^1_{\rm loc}({\Bbb R}_{+})$ and the integral
 \[
\int_{0}^\infty | h(t)| t^{\gamma-1}(1+| \ln t |)^{-\varkappa_{0} } dt< \infty
\]
 converges  for some $\varkappa_{0}\in {\Bbb R} $. In particular, the estimate 
  \[
  | h(t)|\leq C t^{-\gamma } (1+|\ln t|)^{\varkappa}
\]
for some $\varkappa \in {\Bbb R} $ guarantees that $ h\in {\cal S}_{\gamma-1/2}' $. 
  
   The case $\gamma=1$ when $h\in {\cal Z}'_{1/2}$ and hence $\sigma\in {\cal Z}'_{-1/2}$  is most important. It is shown in \cite{Y} that for all bounded Hankel operators $H$, their kernels $h\in {\cal S}'_{1/2}\subset {\cal Z}'_{1/2}$.
The converse is false. For instance, the kernels $h(t)=t^{-1}\ln^k t$ where $k=0, 1,2, \ldots$   satisfy the condition $h\in {\cal S}'_{1/2}$, but the corresponding
Hankel operators   are unbounded if $k\geq  1$ (see \cite{Yf1}).

Note also that the inclusion $h\in {\cal S}_{\gamma-1/2}' $ does {\it not} imply that $\sigma\in {\cal S}_{ 1/2-\gamma}' $. For example,   if $h(t)=  e^{- \alpha t}$, $\Re \alpha >0$  (the Hankel operator $H$  with such kernel has rank $1$), then $h\in {\cal S}'_{\gamma-1/2}$ for all $\gamma>0$, but  if $\Im \alpha\neq 0$ the corresponding function $\sigma(\lambda)= \d (\lambda- \alpha)$ does not belong to ${\cal S}'_{1/2-\gamma}$ for any $\gamma>0$.

The proof of  the main  identity \e{eq:MAIDs} in classes of analytic functions is quite similar to that in Section~2.  First we note an analogue of Proposition~\ref{1Fg}.

 \begin{proposition}\label{1F}
Let $h\in {\cal Z}_{\gamma-1/2}' $ for some $\gamma>0$, and let $\sigma  \in {\cal Z}_{ 1/2-\gamma}'   $ be the corresponding sigma-function. 
    Then      the identity
 \e{eq:MAIDG} 
  holds for arbitrary
   $F\in {\cal Z}_{\gamma-1/2}$.
 \end{proposition}
 
 An analogue of relation \e{eq:gg1} requires a short proof.
     
  \begin{lemma}\label{zz1}
If $f_{1}, f_{2}\in {\cal Z}_{( \gamma -1)/2}$ where $\gamma>0$, then 
   \begin{equation}
{\sf L} (\bar{f}_{1} \star  f_{2} ) = ({\sf L} f_{1})^*{\sf L} f_{2}\in  {\cal Z}_{1/2- \gamma}.
 \label{eq:hh1}\end{equation}
  \end{lemma}

   \begin{pf}
  Since according to \e{eq:heta1}
  \[
|f_{j}(t)| \leq C_{\varkappa} t^{-1+ \gamma/2 } (1+|\ln t|)^{-\varkappa}, \q \forall \varkappa,
 \]
 the integrals $({\sf L} f_{j})(\lambda) $, $j=1,2$,   converge absolutely for all $\lambda>0$. Therefore using the Fubini theorem and making the change of variables $s+t = \tau$, we find that
 \[
   ({\sf L} f_{1} )^*(\lambda)({\sf L} f_{2})(\lambda)= \int_{0}^\infty\int_{0}^\infty e^{-\lambda (t+s)} \bar{f}_{1}(t) f_{2} (s) dtds=\int_{0}^\infty d \tau e^{-\lambda \tau}\int_{0}^\tau \bar{f}_{1}(t ) f_{2} (\tau-t ) d t
     \]
which yields the left-hand side of \e{eq:hh1}. By Corollary~\ref{MEanx}, we have ${\sf L} f_{j} \in {\cal Z}_{(1-\gamma)/2}$, $j=1,2$. Thus the inclusion in  \e{eq:hh1} follows from Lemma~\ref{zz}.
       \end{pf}

       The role of Theorem~\ref{1} is now played by the following result.
       
 \begin{theorem}\label{1z}
Let   $h\in {\cal Z}_{\gamma-1/2}' $  
 for some $\gamma>0$, and let $\sigma \in {\cal Z}_{ 1/2- \gamma}' $ be defined by formula
  \e{eq:LAPL1g}.
  Then       the identity
  \e{eq:MAID} 
  holds for arbitrary
   $f_{1}, f_{2}\in {\cal Z}_{(\gamma-1)/2}$.
 \end{theorem}

 \begin{pf}  
 It suffices to apply identity \e{eq:MAIDG} to $F= \bar{f}_{1}\star f_{2}$ and to use
  Lemma~\ref{zz1}. 
  \end{pf}
  
Thus under the assumption $h\in {\cal Z}_{\gamma-1/2}' $ where $\gamma>0$, the results of Section~2 remain true.

\section{Positive Hankel operators}  

 As we have seen,     the sigma-function $\sigma $  may be a highly singular distribution. However it cannot be too singular for nonnegative Hankel forms. 
   
        \medskip
 
 {\bf 5.1.}
 An important necessary condition of positivity of a Hankel operator $H$ is imposed by 
Bernstein's theorem. Actually, we need its extension to distributions. We consider the problem in a very general setting regarding quadratic forms instead of operators.

\begin{theorem}\label{Bern}
Let $h\in C_{0}^\infty({\Bbb R}_{+})'$ and
\begin{equation}
{\pmb \la}h, \bar{f}\star f {\pmb \ra}\geq 0
\label{eq:BernH}\end{equation}
for all $f\in C_{0}^\infty({\Bbb R}_{+})$. Then there exists a positive measure $M$ on $\Bbb R$ such that
\begin{equation}
h(t)=\int_{-\infty}^\infty e^{-t\lambda} dM(\lambda)  
\label{eq:Bern1}\end{equation}
where  the integral converges for all $t> 0$.
 \end{theorem}
 
 We emphasize that the measure $dM(\lambda) $ may grow almost exponentially as $\lambda\to +\infty$ and it tends to zero super-exponentially as $\lambda\to -\infty$, that is,
 \begin{equation}
 \int_{0}^\infty e^{-t\lambda} dM(\lambda)<\infty \q {\rm and}\q   \int_{0}^\infty e^{t\lambda} dM(-\lambda) <\infty
\label{eq:Bern1x}\end{equation}
  for an arbitrary small $t>0$ and  for an arbitrary large $t>0$, respectively.
 
  Theorem~\ref{Bern} can be viewed as a continuous version of the Hamburger moment problem (see \cite{Hamb} or Theorem~2.1.1 in \cite{AKH}).
  
  Observe that if the function $h(t)$ is a priori supposed to be continuous, then Theorem~\ref{Bern} is exactly the Bernstein theorem on exponentially convex functions (see \cite{Bern} or Theorem~5.5.4 in \cite{AKH}). It is also noted in  \cite{AKH} that due to the theorem of Sierpinski \cite{Sier}, the condition $h\in C ({\Bbb R}_{+})$ in the Bernstein theorem can   be significantly relaxed.
  
   The representation \e{eq:Bern1} is of course  a particular case of \e{eq:conv1}. It is much more precise than \e{eq:conv1} but requires the positivity of ${\pmb\la}h,\bar{f} \star f {\pmb\ra}$. Theorem~\ref{Bern} shows that  the positivity of ${\pmb\la}h,\bar{f} \star f {\pmb\ra}$ imposes very strong conditions on $h(t)$.
  In particular, representation  \e{eq:Bern1} implies that the   distribution $h(t)$ is actually a $C^\infty$ function. It admits the analytic continuation in the half-plane $\Re t>0$ and
\[
h( \tau+ i\sigma)=\int_{-\infty}^\infty e^{-i\sigma\lambda}  e^{-\tau\lambda}dM(\lambda), \q \tau>0.  
\]
This allows us to state the following result.

 \begin{corollary}\label{BernN}
 Under the assumptions of Theorem~\ref{Bern},
the function  $h\in C^\infty ({\Bbb R}_{+})$. Moreover, it  admits the analytic continuation in the right-half plane $\Re t> 0$ and   is uniformly bounded in every strip $\Re t \in (t_{1},t_{2})$ where $0< t_{1} <t_{2} <\infty$.
 \end{corollary}

 Observe that  representation  \e{eq:Bern1} can equivalently be rewritten as
 \begin{equation}
{\pmb \la}h, F {\pmb \ra}=\int_{-\infty}^\infty \ov{({\sf L}F )(\lambda)} dM(\lambda) 
\label{eq:Bern1X}\end{equation}
where the operator ${\sf L}$ is defined by equality \e{eq:LAPj} and $F\in C_{0}^\infty({\Bbb R}_{+})$ is arbitrary. Similarly,
  the  Hankel quadratic form admits the representation
  \begin{equation}
{\pmb \la}h, \bar{f}\star f {\pmb \ra} =\int_{-\infty}^\infty |({\sf L}f)(\lambda)|^2 d M(\lambda), \q \forall  f \in C_{0}^\infty({\Bbb R}_{+}).
\label{eq:QF}\end{equation}
Of course these representations are consistent with formulas \e{eq:MAIDG} and  \e{eq:MAID}.

Our proof of Theorem~\ref{Bern} relies on a reduction to the case of continuous functions $h(t)$. This is   similar in spirit to the extension by L.~Schwartz to distributions of the Bochner theorem on continuous functions of positive type. To be more precise, we follow closely the scheme of \S 3, Chapter~II
of the book \cite{GUEVI}. The difference is that now the Laplace transform plays the role of the Fourier transform and the Laplace convolution defined by \e{eq:HH1} plays the role of the usual convolution. Since  the proof of Theorem~\ref{Bern} is quite far from the mainstream of the present paper, it will be given in the Appendix.

Note that the assertion converse to Theorem~\ref{Bern} is trivially correct: if a function $h(t)$ admits representation \e{eq:Bern1}, then the corresponding Hankel quadratic form is given by relation \e{eq:QF}, and hence it is positive. 


 \medskip
 
 {\bf 5.2.}
 Under assumptions of subs.~4.4 representation \e{eq:Bern1} also holds. In this case one can obtain essentially more detailed information on the measure $dM(\lambda)$.  For the proof of a such result, we  combine Theorem~\ref{1z} with the Bochner-Schwartz theorem (see, e.g., Theorem~3 in \S 3,  Chapter II of the book \cite{GUEVI}). 
 
  It can be stated as follows. Let
a distribution $s\in  {\cal Z}'$ satisfy  the condition
 \begin{equation}
{ \la} s, u^* u {\ra} \geq 0, \q \forall u\in {\cal Z},
\label{eq:HP1}\end{equation}
 (such   $s$ are sometimes called    distributions     of positive type). Then there exists  a nonnegative measure $d{\sf M}(x )$ satisfying the condition
  \begin{equation}
\int_{-\infty}^\infty (1+|x|)^{-\varkappa } d{\sf M}(x )<\infty
\label{eq:Sch}\end{equation}
for some $\varkappa\in {\Bbb R} $ (that is,  of at most polynomial growth at infinity)
and such that   
 \begin{equation}
 { \la} s, \varphi {\ra}=
  \int_{-\infty}^\infty \ov{\varphi(x)}  d{\sf M}(x ), \q \forall \varphi\in {\cal Z}.
\label{eq:HP2}\end{equation}
 In particular, the distribution $s$ can be extended by continuity to the whole Schwartz space $ {\cal S}'$.
 
 Our goal is to prove the following result.

  \begin{theorem}\label{HBx2}
Let $h \in {\cal Z}_{\gamma-1/2}'$ for some $\gamma>0 $ and let condition \e{eq:BernH}
be satisfied  for all
 $f\in {\cal Z}_{(\gamma-1)/2}$. Then the representation
   \begin{equation}
h (t)=\int_{0}^\infty e^{-t\lambda}d M(\lambda), \q \forall t >0, 
\label{eq:C1mm}\end{equation}
holds with
 a   positive measure $dM (\lambda )$ on $  {\Bbb R}_{+}$ satisfying for some $\varkappa \in {\Bbb R}$  the  condition 
      \begin{equation}
\int_{0}^\infty (1+|\ln\lambda|)^{-\varkappa } \lambda^{-\gamma} dM (\lambda )<\infty.
\label{eq:Sch1}\end{equation} 
      \end{theorem} 
      
        \begin{pf}
    Put
           \begin{equation}
 u(x) = e^{-\gamma x/2} ( {\sf L} f) (e^{-x})    
\label{eq:Sch11}\end{equation} 
and
   \[
s(x)= e^{(\gamma -1)x}\sigma (e^{-x})   .
\]
It follows from Corollary~\ref{MEanx} that ${\sf L}  f\in {\cal Z}_{(1-\gamma)/2}$ and hence $u\in {\cal Z}$. Moreover, since
${\sf L} : {\cal Z}_{( \gamma-1)/2} \to {\cal Z}_{(1-\gamma)/2}$  is an isomorphism, for every $u\in {\cal Z}$, we can find $f\in {\cal Z}_{(\gamma-1)/2}$ such that \e{eq:Sch11} holds. According to \e{eq:LL2} we have
 $\sigma = ({\sf L}^*)^{-1}h\in {\cal Z}_{ 1/2-\gamma}'$ and hence $s\in {\cal Z}'$. Making the change of variable $\lambda=e^{-x}$, we see that
\[
 {\pmb\la}\sigma, ({\sf L}  f )^* {\sf L}  f  {\pmb\ra}  = \la s, u^* u\ra.
 \]  
   Therefore using  the main identity  \e{eq:MAID} and assumption \e{eq:BernH}, we obtain  condition \e{eq:HP1}
   on the distribution $s(x)$. 
   
      The Bochner-Schwartz theorem implies that there exists  a positive measure $d{\sf M}(x )$ satisfying  condition \e{eq:Sch} and such that 
representation \e{eq:HP2} holds. Let us now make in \e{eq:HP2} the inverse change of variables $x= - \ln\lambda$ and put
 $\varphi(x)= e^{-\gamma x} \psi (e^{-x})$,
 \[
  e^{-\gamma x}      d{\sf M}(x ) =d M (e^{-x})
  \]
  The measure $d M(\lambda)$ satisfies condition \e{eq:Sch1} and
  \[
  {\pmb \la}\sigma, \psi   {\pmb \ra}=  {  \la}  s, \varphi { \ra}= \int_{0}^\infty \ov{\psi(\lambda)} d M(\lambda).
  \]
  Since $\varphi\in{\cal Z}$ is arbitrary, $\psi\in{\cal Z}_{1/2- \gamma}$ is also arbitrary. Now  the identity 
  \e{eq:MAIDG} with $F={\sf L}^{-1}\psi$ implies the relation 
   \begin{equation}
{\pmb \la}h, F {\pmb \ra}=\int_{0}^\infty \ov{({\sf L}F )(\lambda)} dM(\lambda). 
\label{eq:Bern1Xq}\end{equation} 
    Here $F\in {\cal Z}_{\gamma-1/2} $ is arbitrary because
  ${\sf L} : {\cal Z}_{\gamma-1/2} \to {\cal Z}_{1/ 2-\gamma}$  is an isomorphism. Relations  \e{eq:C1mm} and   \e{eq:Bern1Xq}    are equivalent.
         \end{pf}
         
          \begin{remark}\label{HBb3}
       If $h \in {\cal Z}_{\gamma_{1}-1/2}' \cap {\cal Z}_{\gamma_{2}-1/2}'$ for some $0<\gamma_{1} <\gamma_{2} <\infty $, then  the representation
\e{eq:C1mm} 
holds with
 a     measure $dM (\lambda )$  satisfying instead of \e{eq:Sch1} the stronger condition
       \[
\int_{1}^\infty (1+|\ln\lambda|)^{-\varkappa_{1} } \lambda^{-\gamma_{1}} dM (\lambda )+ \int_{0}^1 (1+|\ln\lambda|)^{-\varkappa_{2} } \lambda^{-\gamma_{2}} dM (\lambda )<\infty
\]
 for some $\varkappa_{1}, \varkappa_{1} \in {\Bbb R}$.
             \end{remark}

    Under the assumptions of Theorem~\ref{HBx2}, $h(t)$ satisfies the conclusions of Corollary~\ref{BernN}. Furthermore, we have
      
      \begin{corollary}\label{HBx3}
      Under the assumptions of Theorem~\ref{HBx2}, for all $t>0$ and all $n=0,1,2,\ldots$,  inequalities
    \begin{equation}
(-1)^n h^{(n)}(t)\geq 0
\label{eq:CoMo}\end{equation}
hold    $($such functions $h(t)$ are   called completely monotonic$)$.  Moreover, for some $\varkappa\in {\Bbb R}$ and $C>0$ we have the estimate  
    \begin{equation}
 h(t) \leq C t^{-\gamma} (1+|\ln t|)^\varkappa ,\q t>0.
\label{eq:CoMo1}\end{equation}
      \end{corollary} 
      
      All these assertions are direct consequences of the representation \e{eq:C1mm}. In particular, under condition  \e{eq:Sch1} we have
      \[
      h(t) \leq C \max_{\lambda\geq 0} \big(e^{-t\lambda}\lambda^\gamma (1+|\ln \lambda|)^\varkappa\big)
      \]
      which yields \e{eq:CoMo1}.

          Under the assumptions of Theorem~\ref{HBx2},
 we have the representation
    \[
{\pmb\la}h,\bar{f}  \star f {\pmb\ra} =  \int_0^\infty |({\sf L} f ) (\lambda) |^2 d M(\lambda), \q \forall f \in {\cal Z}_{(\gamma-1)/2}.
\]
  In contrast to  \e{eq:QF} the integral here is taken over the positive half-line only. This is of course due to stronger assumptions on $h(t)$.

 Note that according to the  Bernstein theorem (see the original paper  \cite{Bern} or Theorems~5.5.1 and 5.5.2 in the book \cite{AKH}) condition \e{eq:CoMo} implies that the function $h(t)$ admits the representation \e{eq:C1mm} with some measure $dM(\lambda)$. Of course, condition \e{eq:CoMo}  does not impose any restrictions on the measure $dM(\lambda)$ (except that the integral \e{eq:C1mm} is convergent for all $t>0$). In contrast to the  Bernstein theorem we deduce the representation \e{eq:C1mm} from the positivity of the Hankel form. In this context  condition  \e{eq:Sch1} is due to the assumption $h \in {\cal Z}'_{\gamma-1/2}$.

   We also  mention  that  H.~Widom  considered in \cite{Widom} Hankel operators $H$ with kernels $h(t)$ admitting the representation \e{eq:C1mm}. He showed that $H$ is bounded if and only if $M([0,\lambda))=O(\lambda)$ as $\lambda\to 0$ and as $\lambda\to \infty$.  In this case $h(t)\leq C t^{-1}$ for some $C>0$. To a certain extent, Theorem~\ref{HBx2}
   and estimate  \e{eq:CoMo1} 
  can be regarded as an extension of Widom's results to unbounded operators.

\section{Quasi-Carleman operators}  

{\bf 6.1.}
Here we consider Hankel operators $H$  (we call them ``quasi-Carleman"   operators) with kernels \e{eq:E1r}  that belong to the set $ C_{0}^\infty ({\Bbb R}_{+})'$  for all $\alpha\in {\Bbb R}$, $r\geq 0$ and $k\in{\Bbb R}$. To be more precise, we study the corresponding quadratic forms.
It can be shown that these forms give rise to self-adjoint operators $H$  in the space  $ L^2 ({\Bbb R}_{+})$ if and only if $\alpha >0$ or $\alpha =0$, $k>0$ (in these cases $h(t)\to 0$ as $t\to \infty$).
Moreover, if $\alpha>0$, $r>0$, then $H$ is compact for all $k$. If $\alpha>0$, $r=0$, then it  is compact for   $k>-1$ and is bounded for   $k= -1$. If $\alpha=0$, $r>0$, then it  is compact for   $k<-1$ and is bounded for   $k= -1$. Finally, if
$\alpha= r=0$, then $H$ is bounded if and only if $k=-1$. In all these cases we   have the equality $N_{\pm}(H)= N_{\pm}(h)$.

There are probably no chances to explicitly find the spectrum and eigenfunctions of quasi-Carleman operators. The only exceptions are the cases  $k=-1$, $\alpha=0$ (if in addition $r= 0$, then $H$ 
 is the Carleman operator) and $k=-1$, $r=0$, $\alpha > 0$ considered by F.~G.~Mehler \cite{Me} and W.~Magnus \cite{Ma}, respectively (see also  \S 3.14 of the book \cite{BE} and the papers \cite{Ro}, \cite{Y1}).

Our first goal is to prove formula \e{eq:bbr7} for the sigma-functions. We consider all  $k\in {\Bbb C}$ and start with the case $\Re k<0$ when distribution \e{eq:bbr7}   does  not have a strong singularity at the point $\lambda=\alpha$.
Formally, the proof is quite simple. Indeed, for $h(t)= t^k$, we apply  the relation
\[
\int_{0}^\infty \lambda^{-k-1} e^{-\lambda t} d\lambda = \Gamma (-k) \, t^{k} .
\]
To pass to the general case, one can use  the following observation. If
\[
h_{r,\alpha}(t) =h(t+r)e^{-\alpha t}, \q r\geq 0,  
\]
(that is, a kernel $h(t)$ is shifted and multiplied by an   exponential), 
then according to  \e{eq:conv1}  the corresponding sigma-function equals
\[
  \sigma_{r,\alpha}(\lambda)=  e^{-r (\lambda-\alpha)} \sigma (\lambda-\alpha).
\]


Let us now give the precise proof of \e{eq:bbr7}.

    \begin{lemma}\label{MIS}
If $\alpha\in {\Bbb R}$, $r\geq 0$ and $\Re k<0$, then   for all $F\in C_{0}^\infty ({\Bbb R}_{+})$ the identity 
 \begin{equation}
 \int_0^\infty  (t+r)^{k } e^{-\alpha t} \overline{F(t)}  dt =
  \frac{e^{\alpha r}}{\Gamma( -k)} \int_{\alpha}^\infty (\lambda-\alpha)^{-k-1 }e^{-r \lambda} \ov{({\sf L}F)(\lambda) } d \lambda
\label{eq:C1}\end{equation}
 holds.
 \end{lemma}
 
 \begin{pf}
 We use   definition \e{eq:LAPj} of the operator ${\sf L}$ and according to the Fubini theorem interchange the order of integrations in  the right-hand side of \e{eq:C1}.  Thus it  equals
  \[
 \frac{e^{\alpha r}}{\Gamma( -k)}    \int_0^\infty dt    \overline{F(t)}    
  \int_{\alpha}^\infty (\lambda-\alpha)^{-k-1 }e^{-(t+r) \lambda}d\lambda.
\]
Since the integral over $\lambda$ equals $\Gamma( -k ) (t+r)^{k} e^{-\alpha ( t+r)}$, this yields  the left-hand side of \e{eq:C1}.  
\end{pf}

Our next goal is to extend formula \e{eq:C1} to $k$ in the right-half plane. The left-hand side of \e{eq:C1} is obviously an analytic function of  $  k \in {\Bbb C}$. As is well known, the analytic continuation of the integral in  the right-hand side of \e{eq:C1} to the strip 
  $n<\Re k  <n+1$ where $ n\in{\Bbb Z}_{+}$ is given by the integral
  \begin{equation}
  \int_{\alpha}^\infty (\lambda-\alpha)^{-k-1 }  \big( \omega(\lambda) -\sum_{p=0}^{n}  \frac{1}{p!} \omega^{(p) }(\alpha) (\lambda-\alpha)^p \big) d\lambda=:   \int_0^\infty (\lambda-\alpha)_{+}^{-k-1 }   \omega(\lambda)   d\lambda
\label{eq:di}\end{equation}
where $\omega(\lambda)= e^{-r \lambda} \ov{({\sf L}F)(\lambda)} $. Here we use the standard notation $(\lambda-\alpha)_{+}^{-k-1 } $ for the distribution determined by this formula (we refer, for example, to the book \cite{GUECH} for a discussion of such distributions). This distribution is also  well defined, although by a slightly different formula,  on the lines $\Re k\in {\Bbb Z}_{+}$.

This concludes the proof of relation \e{eq:bbr7}. Let us formulate the result obtained.

\begin{lemma}\label{Hks}
Let $h(t)$ be given by formula \e{eq:E1r} where  $\alpha\in {\Bbb R}$, $r\geq 0$. If $k \in{\Bbb R}\setminus{\Bbb Z}_{+}$, then the sigma-function is given by equality \e{eq:bbr7}. If $k \in{\Bbb Z}_{+}$, it is given by equality \e{eq:dii}. 
    \end{lemma}
    
    Putting together this result with Theorem~\ref{1}, we get the following assertion.
    
    \begin{proposition}\label{HM}
Let $h(t)$ be given by formula \e{eq:E1r}, and let the function $\sigma(\lambda)$ be given by equalities \e{eq:bbr7} or \e{eq:dii}. Then the identity  \e{eq:MAID} holds for all $f_1, f_{2}\in C_{0}^\infty ({\Bbb R}_{+})$.
    \end{proposition}

Since for the sigma-function \e{eq:bbr7}, the function
 \[
  s (x) := \sigma (e^{-x})=
   \frac{e^{\alpha r}}{\Gamma(-k)}   (e^{-x} -\alpha)_{+}^{-k-1}e^{-r e^{-x}} 
\]
belongs to the Schwarz class  $  {\cal S}'$, we see that $\sigma\in {\cal S}'_0$.
The same is true for the sigma-function \e{eq:dii}.

\medskip
 
 {\bf 6.2.}  
  Our next goal is to calculate the numbers $N_\pm ( h):= N_\pm ( h; C_{0}^\infty ({\Bbb R}_{+} ))$. Observe that the number $N_\pm ( h  )$ does not depend on $\alpha\in {\Bbb R} $ in definition \e{eq:E1r}. Indeed, if $\pm {\pmb \la }h, \bar{f}\star f{\pmb \ra } >0$ for all $f\neq 0$ in some linear space ${\sf D}\subset C_{0}^\infty ({\Bbb R}_{+} )$ and $h_{\gamma}(t)= e^{\gamma t} h(t)$ for some
  $\gamma\in {\Bbb R} $, then  $\pm {\pmb \la }h_{\gamma }, \bar{f}\star f{\pmb \ra } >0$  for all $f\neq 0$ in the  linear space ${\sf D}_{\gamma}$ consisting of functions $e^{-\gamma t} f(t)$ where $f\in{\sf D}$. The  spaces ${\sf D} $  and ${\sf D}_{\gamma}$ have of course the same dimension.
  
The case  of Hankel operators of finite rank was treated  in \cite{Yf}. If $k\in  {\Bbb Z}_{+}$ and $\alpha>0$, then the form $  {\pmb \la }h, \bar{f}\star f{\pmb \ra } $ gives  rise to a Hankel operator of  rank $k+1$ and $N_{+} (H)= N_\pm ( h )$. The above remark allows us to extend the result of  \cite{Yf} to all 
  $\alpha\in {\Bbb R} $. Let us state the corresponding assertion.

  \begin{theorem}\label{HKFR}
Let $h(t)$ be given by formula \e{eq:E1r} where   $\alpha\in {\Bbb R} $, $r\geq 0$ and  $k  \in {\Bbb Z}_{+}$. 
Then  $N_+ ( h )= N_-( h ) +1= k/2
+1$ if $k$ is even and
  $N_\pm ( h ) =(k+1)/2$ if $k$ is odd.
   \end{theorem}
   
   Our goal here is to prove the following result.
   
    \begin{theorem}\label{HKL}
Let $h(t)$ be given by formula \e{eq:E1r} where   $\alpha \in {\Bbb R}$ and $r\geq 0$. If $k<0$, then  ${\pmb \la }h, \bar{f}\star f{\pmb \ra }\geq 0$ for all $f\in C_{0}^\infty ({\Bbb R}_{+} )$. If
  $k >0$ but $k \not \in {\Bbb Z}_{+}$, then  $N_+ ( h )=  [k] /2+1$, $N_- ( h )= \infty$   for even $[k]$ and  $N_- ( h )= ([k]+ 1)/2$, $N_+ ( h )= \infty$ for odd $[k]$. 
 \end{theorem}
  
     If $k<0$, then, by formula  \e{eq:bbr7}, $\sigma\in L^1_{\rm loc}$ and   $\sigma(\lambda)\geq 0$. Therefore it suffices to use the identity \e{eq:MAID}.

The   case $k > 0$  is essentially more complicated. Without loss of generality, we   suppose that $\alpha>0$. Then the operators $H$ are compact and  $N_\pm (H)=N_\pm (h )$. We proceed from the  assertion which follows from  Theorem~\ref{HByz} if    formula  \e{eq:bbr7} for $\sigma (\lambda)$ is taken into account.

  \begin{lemma}\label{Hkl}
  Let $h(t)$ be given by formula \e{eq:E1r} where   $\alpha >0$,  $r\geq 0$  and $k \in {\Bbb R}_{+} \setminus {\Bbb Z}_{+} $.
Define the form   
 \begin{equation}
\sigma[ w,w]=    \frac{e^{\alpha r}}{\Gamma(-k)}   \int_{0}^\infty (\lambda-\alpha)_{+}^{-k-1} e^{-r\lambda} |w(\lambda) |^2 d\lambda 
\label{eq:L4}\end{equation}
on functions $w \in C_0^\infty ({\Bbb R}_{+})$ and put $N_\pm (\sigma ) : = N_\pm (\sigma;  C_0^\infty ({\Bbb R}_{+}))$.  Then $N_\pm (h ) = N_\pm (\sigma )$.
 \end{lemma}
 
 Below we need two   elementary assertions on distributions $\mu^{-k-1}_+ $.
 
  \begin{lemma}\label{SiDx}
     Let $\alpha>0$, $k \in {\Bbb R}_{+} \setminus {\Bbb Z}_{+} $ and $n=[k]$.
Suppose that a function $ w \in C_{0}^\infty({\Bbb R}_{+})$ and
 \begin{equation}
w (\alpha) =w' (\alpha)=\cdots =w^{(\ell -1)} (\alpha) =0,
\label{eq:L3f}\end{equation}
where 
 \begin{equation}
\ell= \ell (n) =  n/2+1\; {\rm for} \; n  \; {\rm even}\q {\rm and } \q\ell= \ell (n) =  (n+1)/2 \; {\rm for} \; n \, {\rm odd}.
\label{eq:L3v}\end{equation}
Then form \e{eq:L4} satisfies the inequality
 \begin{equation}
(-1)^{n+1}\sigma [w, w] \geq 0.
\label{eq:L3fj}\end{equation}
 \end{lemma}
 
 \begin{pf}
 Put $\varphi (\lambda)  = e^{-r\lambda} |w(\lambda) |^2$. It follows from
   definition \e{eq:di} that
   \begin{equation}
\Gamma(-k) e^{- \alpha r} \sigma[ w,w]
  =  \int_{\alpha}^\infty (\lambda-\alpha)^{-k-1} \big(\varphi (\lambda)  
  - \sum_{p=0}^n  \frac{1}{p!   }  \varphi^{(p)} (\alpha) (\lambda-\alpha)^p\big)d\lambda.
\label{eq:L3}\end{equation}
Under assumptions \e{eq:L3f},  \e{eq:L3v}  we have
\begin{equation}
\varphi (\alpha) =\varphi' (\alpha)=\cdots =\varphi^{(n)} (\alpha) =0
\label{eq:E1Gr}\end{equation}
so that  the right-hand side of \e{eq:L3} is   nonnegative.
 Since $\Gamma(-k)<0$ for $n$  even   and $\Gamma(-k)> 0$ for  $n$ odd, equality \e{eq:L3} implies  \e{eq:L3fj}.
  \end{pf}

     \begin{lemma}\label{SiD}
     Suppose that $\alpha>0$ and $k \in {\Bbb R}_{+} \setminus {\Bbb Z}_{+} $.
Let a   function $\psi=\bar{\psi} \in C_{0}^\infty ({\Bbb R}_{+})$  satisfy the conditions
 \begin{equation}
\psi (\alpha) = 1 \q  \mbox{and}\q \psi' (\alpha) = \cdots =\psi^{(n)} (\alpha) = 0  \q  \mbox{if} \q n =[k]\geq 1,
\label{eq:L6}\end{equation}
and let $Q(\lambda)$ be a polynomial of $\deg Q \leq n$. Then
 \begin{equation}
   \int_{0}^\infty (\lambda-\alpha)^{-k-1}_+ Q(\lambda)\psi^2 (\lambda) d \lambda
   =   \int_{\alpha}^\infty (\lambda-\alpha)^{-k-1} Q(\lambda) ( \psi^2 (\lambda) -1) d\lambda. 
   \label{eq:did}\end{equation}
 \end{lemma}
 
 \begin{pf}
 Put $\omega (\lambda) =Q(\lambda ) \psi^2 (\lambda)$.
 According to \e{eq:L6}   we have
$\omega^{(p)} (\alpha) =  Q^{(p)} (\alpha)$ for all $ p=0,\ldots, n$,
whence
\[
 \sum_{p=0}^n  \frac{1}{p!   } \omega^{(p)}  (\alpha) (\lambda-\alpha)^p = 
 \sum_{p=0}^n  \frac{1}{p!   } Q^{(p)} (\alpha) (\lambda-\alpha)^p =  Q(\lambda)
\]
 if $n\geq \deg Q$. Therefore relation \e{eq:did} is a direct consequence of definition  \e{eq:di}.
   \end{pf}

 Now we are in a position to prove Theorem~\ref{HKL}. 
 In view of Lemma~\ref{Hkl}, to that end we   only have
   to calculate the numbers $ N_\pm (\sigma)$. Let us  consider $ N_+ (\sigma)$ for even $n$ and  $ N_- (\sigma)$ for odd $n$. First we show that $N_\pm (\sigma) \leq \ell$ with $\ell$ defined by \e{eq:L3v}. Suppose the contrary.  Then  there exist  linearly independent functions $ w_j\in C_{0}^\infty({\Bbb R}_{+})$, $j=1,\ldots,\ell +1$,  such that  
      \begin{equation}
 (-1)^{n+1} \sigma [w,w] <0
  \label{eq:M3ff}\end{equation}
  on all their nontrivial linear combinations  
   \begin{equation}
  w (\lambda) =\sum_{j=1}^{\ell +1} c_j w_j (\lambda).
  \label{eq:L3ff}\end{equation}
 Substituting this expression into $\ell$ equations \e{eq:L3f}, we   find a nontrivial solution of this system for the coefficients $c_{1},\ldots,  c_{\ell +1}$.  According to Lemma~\ref{SiDx} for the corresponding function \e{eq:L3ff} we have inequality \e{eq:L3fj}. 
  Clearly, $w\neq 0$ because the  functions $w_{1},\ldots, w_{\ell +1}$ are linearly independent.
  Therefore inequalities \e{eq:L3fj} and \e{eq:M3ff} are incompatible.
 
Let us prove the opposite estimates $N_\pm (\sigma) \geq \ell$. We choose a function  $\psi = \bar{\psi} \in C_0^\infty  ({\Bbb R}_{+})$ satisfying conditions \e{eq:L6} and such that $0\leq \psi(\lambda)\leq 1 $.
Let us calculate form  \e{eq:L4} on functions 
  \begin{equation}
w (\lambda) =   P (\lambda ) \psi  (\lambda  ) e^{r\lambda/2}
\label{eq:L5}\end{equation}
where $P (\lambda) $ is an arbitrary polynomial of $\deg P  \leq [ n/2 ]$.
Applying   Lemma~\ref{SiD} to $Q(\lambda)=|P(\lambda)|^2$,  we see 
  that
  $\Gamma(-k)e^{-\alpha r} \sigma [w,w]$ equals expression \e{eq:did}.
 This yields a linear subspace of functions \e{eq:L5} of dimension $[n/2] +1$ where   $\Gamma(-k)  \sigma [w,w]<0$ for all $P\neq 0$.

Thus we have proven that $N_+ (h)=N_+ (\sigma) =n/2 + 1$  for $n$  even  and  $N_- (h)=N_- (\sigma) = (n+1)/2$ for $n$  odd. Since $\Gamma (-k) \sigma (\lambda)>0$ for all $\lambda>\alpha$, it follows from part $2^0$ of Theorem~\ref{HBy} that $N_- (h)= \infty$  for   $n$ even and  $N_+(h) = \infty$ for   $n$ odd (this result also follows from the fact that the rank of $H$ is infinite). 
The proof of Theorem~\ref{HKL} is complete.

  \medskip
 
 {\bf 6.3.}  
 The proof of Theorem~\ref{HKL} actually  relies  only on the study of the singularity of the sigma-function at the point $\alpha>0$.  To emphasize this idea, we obtain here  more general results where  conditions are formulated in terms of the sigma-function $\sigma(\lambda)$ of Hankel operators without making specific assumptions on their kernels $h(t)$. To obtain an upper bound on numbers \e{eq:abc}, we require that the singularity of $\sigma(\lambda)$ at $\lambda=\alpha$ is not too strong. 
 
  \begin{lemma}\label{CB}
Suppose that $h\in C_{0}^\infty ({\Bbb R}_{+})'$ and that the  corresponding sigma-function $ {\sigma}  (\lambda)$ is continuous away from the point $\lambda=\alpha$, bounded as $\lambda\to 0$ and $\lambda\to \infty$ and, for some $n\in {\Bbb Z}_{+}$,  the function $(\lambda-\alpha)^{n+1}  {\sigma} (\lambda)$ belongs to $  L^1_{\rm loc} ({\Bbb R}_{+})$.
Assume also that 
  \begin{equation}
(- 1 )^{n+1}  \sigma  (\lambda)  \geq 0,\q \lambda\neq\alpha .
\label{eq:E1rs11}\end{equation}
Then $N_+ (h )\leq n/2+1$  for
 $n$   even and  $N_- (h )\leq (n+ 1)/2$  for $n$  odd.
\end{lemma}

\begin{pf}
  If a function $w\in C_{0}^\infty ({\Bbb R}_{+} )$   satisfies conditions \e{eq:L3f} where $\ell$ is defined by \e{eq:L3v}, then the function $\varphi(\lambda) = |w (\lambda )|^2$   satisfies conditions \e{eq:E1Gr} so that $\varphi(\lambda)= O(|\lambda-\alpha|^{n+1} )$.
It follows that
 \[
(-1)^{n+1} \sigma [ w, w]=(-1)^{n+1}  \int_{0}^\infty\sigma   (\lambda)   | w (\lambda)|^2d\lambda 
\]
where the integral converges (at the point $\lambda=\alpha$).  By condition \e{eq:E1rs11} this expression   is nonnegative. So it remains to repeat the proof of Theorem ~\ref{HKL} of the upper bounds on the numbers $N_{\pm} (\sigma)$.
\end{pf}
 
 To obtain a lower bound on the numbers $N_\pm (h ) $,
 we assume that
 \begin{equation}
\sigma (\lambda)= \sigma_0 (\lambda)+ \ti{\sigma} (\lambda) 
\label{eq:E1rs1}\end{equation}
 where $\sigma_0(\lambda)$  is given by formula \e{eq:bbr7} and  the singularity of $\ti{\sigma} (\lambda)$ at the point $\alpha$ is weaker than that of $\sigma_0(\lambda)$.
  Namely, we accept the following
 
 \begin{assumption}\label{HG}
 Set $\varphi_\varepsilon(\lambda)= \omega ((\lambda-\alpha)/ \varepsilon) $ where $\omega\in C_0^\infty ({\Bbb R})$. Then
\[
\la \ti{\sigma}, \varphi_\varepsilon\ra =o (\varepsilon^{-k}) \q {\rm as}\q \varepsilon\to 0.
\]
 Moreover, it is supposed   that this relation holds  uniformly for functions $\omega$  having common support in ${\Bbb R} $ and uniformly bounded in $C^n$-norm (as usual $n=[k]$).
 \end{assumption}

    The following result generalizes Theorem~\ref{HKL}.

\begin{theorem}\label{HKLG}
In addition to the conditions of Lemma~\ref{CB}, suppose that representation \e{eq:E1rs1} holds with $\sigma_0 $  given by formula \e{eq:bbr7} where $k\in (n,n+1)$ and $\ti{\sigma}$ satisfying Assumption~\ref{HG}. Then $N_+ (h ) = n/2+1$  for
 $n$   even and  $N_- (h ) = (n+ 1)/2$  for $n$  odd. 
 \end{theorem}

\begin{pf}
The upper estimate on the numbers $N_{\pm} (h)$ is given by Lemma~\ref{CB}.

 To prove the lower estimate, we use again test functions \e{eq:L5} but introducing a small parameter $\varepsilon $ we put
  \begin{equation}
w_{\varepsilon} (\lambda) =   P ((\lambda-\alpha)/\varepsilon)\psi (\alpha+(\lambda-\alpha)/\varepsilon) e^{r\lambda/2}.
\label{eq:L5e}\end{equation}
Here $\psi =\bar{\psi}\in C_{0}^\infty ({\Bbb R}_{+}) $ satisfies conditions \e{eq:L6}, $0\leq \psi(\mu)\leq 1$ and $P(\mu)$ is an arbitrary polynomial of $\deg P\leq [n/2]$.  Similarly to the proof of \e{eq:did}, we now find that
\[
\Gamma(-k)e^{-\alpha r} \sigma_{0} [w_{\varepsilon},w_{\varepsilon}]=  \varepsilon^{-k}    \int_{0}^\infty  \mu^{-k-1}  | P ( \mu )|^2
  \big(\psi(\alpha+\mu)^2- 1 \big)   d \mu.
\]
Put $\| P\|^2 = |p_0|^2+\cdots + |p_{ [n/2]}|^2 $ where $p_{0}, p_{1}, \ldots, p_{ [n/2]}$ are the coefficients of $P(\mu)$. Since
 \[
 \max_{ \|P\| =1}   \int_{0}^\infty  \mu^{-k-1}  | P (\mu )|^2
  \big(\psi(\alpha+ \mu)^2- 1 \big)    d\mu  <0,
 \]
 we see that
\begin{equation}
\Gamma(-k)  \sigma_{0} [w_{\varepsilon},w_{\varepsilon}]\leq   -c \varepsilon^{-k}    \|P \|^2 
\label{eq:L7e1}\end{equation}
 for some constant $c>0$.
Applying  Assumption~\ref{HG} to the function $\omega(\mu)= | P (\mu ) |^2 \psi^2(\alpha+\mu) e^{  r\mu \varepsilon}$, we see that
 \begin{equation}
 \ti{\sigma} [ w_\varepsilon,  w_\varepsilon] = o( \varepsilon^{-k}   )
\label{eq:L7e2}\end{equation}
where the limit is uniform for all polynomials with $\|P \| \leq 1$.
Combining estimates \e{eq:L7e1} and \e{eq:L7e2}, we see that
  \[
\Gamma(-k)  \sigma [ w_\varepsilon,  w_\varepsilon]  \leq -c \|P \|^2  \varepsilon^{-k} (1+ o ( 1)), \q c>0  .
\]
Since the right-hand side here is negative for sufficiently small $\varepsilon$, this yields us a space of dimension $ [n/2]+1$ where the form $\Gamma(-k) \sigma$ is negative.
  \end{pf}

 
 We note that if the  function   $\sigma(\lambda)$  changes the sign for $\lambda\neq\alpha$,  then   $N_\pm (h )= \infty$ according to Theorem~\ref{HBy}. 
 

 \begin{example}\label{HKLG1}
 Let
 \begin{equation}
h(t)=\big(  ( t+r)^{k}+ \sum_{j=1}^{j_0} a_j ( t+r)^{k_j}\big) e^{-\alpha t} ,\q a_j =\bar{a}_j,   \q r\geq 0,
\label{eq:E1rs}\end{equation}
where $k_{j}\in [0,k)$ for all $ j=1,\ldots, k_{j_0}$.  According to
 formula \e{eq:bbr7}  representation \e{eq:E1rs1} is satisfied with
  \begin{equation}
 \ti{\sigma} (\lambda) = e^{-r (\lambda-\alpha)}\sum_{j=1}^{j_0} \frac{a_j}{\Gamma(-k_j)} (\lambda-\alpha)_{+}^{-k_j }
\label{eq:E2rs}\end{equation}
(if $k_{j}\not\in {\Bbb Z}_+$ for all $ j=1,\ldots, k_{j_0}$).  
 Assumption~\ref{HG} holds true  because all $ k_{j} <k  $. Now condition \e{eq:E1rs11}  is fulfilled if
   \begin{equation}
1 + \sum_{j=1}^{j_0} a_j \frac{\Gamma(-k) }{\Gamma(-k_j)}\mu^{k-k_j } \geq 0, \q \forall \mu>0.
\label{eq:E2rsx}\end{equation}
 In particular, it suffices to require that  $(-1)^{n_{j}-n} a_{j}\geq 0$ where $n_{j}=[k_{j}]$ for all $j$.   Then all conclusions of Theorem~\ref{HKLG}  hold. 
  \end{example}
 

We emphasize that it is allowed in \e{eq:E1rs} that $k_j \in {\Bbb Z}_+$. According to Lemma~\ref{Hks} the sigma-functions $\sigma_{j}(\lambda)$ of such kernels $t^{k_j} e^{-\alpha t} $ are combinations of delta functions $\d(\lambda-\alpha)$  and  their derivatives   so that $\sigma_{j}(\lambda)=0$  for $\lambda\neq\alpha$.  Therefore the corresponding term in  \e{eq:E2rsx} should be omitted.

\medskip
 
 {\bf 6.4.}  
Finally, we consider the case when the sigma-function has singularities at several points. It turns out that  the contributions of different  singularities   to the numbers $N_\pm (h)$ are independent of each other. 

\begin{theorem}\label{SNDb}
Let
  \begin{equation}
h(t)= \sum_{m=1}^M   (-1)^{n_m+1} b_m (t+r_m)^{k_m} e^{-\alpha_m t} ,   \q b_m= \bar{b}_m, \; r_{m}\geq 0,
\label{eq:E1G1}\end{equation}
where $k_{m}> 0$, $k_{m}\not\in{\Bbb Z}_{+}$ and $n_m=[k_{m}]$. Then:

$1^0$ If $b_{m}>0$ for some $m=1,\ldots, M$, then $N_{+} (h)=\infty$. If $b_{m}<0$ for some $m=1,\ldots, M$, then $N_{-} (h)=\infty$. 

$2^0$ Put
 \begin{equation}
{\cal N}= \sum_{m=1}^M \, [n_{m}/2] + M .
\label{eq:SBR}\end{equation}  
 If all $b_{m}<0$, then $N_{+}(h)= {\cal N}$.  If all $b_{m}>0$, then $N_{-}(h)= {\cal N}$. 
  \end{theorem} 
  
  \begin{pf}
  It follows from formula \e{eq:bbr7} that the sigma-function of kernel \e{eq:E1G1} equals
   \[
   \sigma(\lambda) =\sum_{m=1}^M  \sigma_{m}(\lambda)\q {\rm where} \q
 \sigma_{m}(\lambda)=  \frac{(-1)^{n_m+1} b_m}{\Gamma(-k_{m})} (\lambda-\alpha_{m})_{+}^{-k_{m}-1} e^{-r_{m} (\lambda-\alpha_{m})}.
\]
Obviously, the function  $   \sigma(\lambda)$ is continuous away from the points $\alpha_{1},\ldots, \alpha_{M}$.  Clearly, $\sigma\in {\cal S}_{0}$ if $\alpha_{m}> 0$ for all $m=1,\ldots, M$ which we can always suppose.
 Note that $(-1)^{n_m+1}  \Gamma(-k_{m})>0$. Therefore if $b_{m}>0$ ($b_{m}<0$) for some $m$, then this function tends to $+\infty$ ($-\infty$) as $\lambda\to \alpha_{m}+0$. In this case, by Theorem~\ref{HBy}, the positive (negative) spectrum of the operator $H$ is infinite.

Let us prove part  $2^0$. Suppose, for example, that $b_{m}>0$ for all $m=1,\ldots, M$; then $ \sigma (\lambda)\geq 0$ for all $\lambda\not\in\{\alpha_{1}, \ldots, \alpha_{m}\}$. Let a function $w\in C_{0}^\infty ({\Bbb R}_{+})$   satisfy the conditions $w^{(p)}(\alpha_{m})=0$ for $p=0,1, \ldots, [n_{m}/2] $    and all $m=1,\ldots, M$. Then according to Lemma~\ref{SiDx} we have $\sigma_{m}[w,w]\geq  0$ 
for all $m=1,\ldots, M$ and hence $\sigma [w,w]\geq  0$. Quite similarly to the proof of the upper bound on  $N_{-}(h)$ in Theorem~\ref{HKL}, this implies that $N_{-}(h)\leq {\cal N}$. 

To prove the opposite inequality, we consider trial functions $w_{\varepsilon,m}$, $m=1,\ldots, M$, defined by formula \e{eq:L5e} where $\alpha=\alpha_{m}$, $r=r_{m}$ and $P_{m}$ is a polynomial of degree $[ n_{m} /2]$. Instead of \e{eq:L7e1}, we now have the estimate
\[
  \sigma_m [w_{\varepsilon,m},w_{\varepsilon,m}]\leq   -c \varepsilon^{-k_{m}}    \|P_{m} \|^2 ,\q c>0.
\]
Since the functions $ \sigma_p (\lambda)$ where $p\neq m$ are continuous at the point $\alpha_{m}$, this implies the same estimate on $  \sigma [w_{\varepsilon,m},w_{\varepsilon,m}]$. Taking linear combinations of functions  $w_{\varepsilon,m}$, we obtain a space of dimension $\cal N$ where the form $\sigma$ is negative.
      \end{pf}

      \begin{remark}\label{SNDbx}
Suppose that some of $k_{m}$ in  
\e{eq:E1G1} are integers. Then $N_{\pm} (h)=\infty$ if $\pm b_{m}>0$ for some $m $ such that $k_{m} \not\in{\Bbb Z}_{+}$. Assertion $2^0$ of Theorem~\ref{SNDb}  remains unchanged.
  \end{remark} 
  
    \begin{remark}\label{SNDbK}
Let $h(t)$ be given by formula \e{eq:E1G1}. Then $\pm {\pmb\la}h, \bar{f} \star f {\pmb\ra}\geq 0$ for all $f\in C_{0}^\infty ({\Bbb R}_{+})$ if and only if $k_{m}\leq 0$ and $\mp (-1)^{n_{m}}b_{m} \geq 0$ for all $m=1,\ldots, M$.
  \end{remark}

      \section{Discrete representation }
      
      
   Here we consider Hankel operators $Q$  defined      by   equality  \e{eq:K1} in the space     $ l^2({\Bbb Z}_{+})$ of sequences $g =(g_{0}, g_{1}, \ldots)$ and discuss their relation by formula \e{eq:K2}   to integral Hankel operators $H$    in the space     $ L^2({\Bbb R}_{+})$.
   It turns out that the concept of the sigma-function is also very convenient for finding  a link between matrix elements $q_{n}$ of $Q$ and kernels $h(t)$ of $H$.
   
      
      \medskip
      
{\bf 7.1}.
Similarly to the continuous case, the most general definition of Hankel operators $Q$ is given  in terms of quadratic forms 
\[
q[g,g]=  \sum_{n,m=0}^\infty q_{n+m} g_{m}\bar{g}_{n} 
\]
 considered on the set $\ell_{0}\subset l^2({\Bbb Z}_{+})$ of elements $g$ with only finite number of non-zero components. This definition does  not require any assumptions on elements $q_{n}$, but it does not guarantee that $Q$ is correctly defined as an  operator (even unbounded)  in $l^2 ({\Bbb Z}_{+})$. 
 
 Let us construct  a unitary operator ${\bf U} : l^2({\Bbb Z}_{+})\to  L^2({\Bbb R}_{+})$ such that relation \e{eq:K2} is satisfied. To be precise, we consider
  the relation 
\begin{equation}
{\pmb \la} h,\ov{ {\bf U}  g }\star {\bf U} g  {\pmb \ra}=  q[g,g], \q g\in \ell_{0},
\label{eq:KH}\end{equation}
between the corresponding quadratic forms.

 Recall that, for an arbitrary value  of the parameter $\kappa >-1$,   the Laguerre polynomial (see \cite{BE}, Chapter~10.12) of degree $n$ is defined by the formula
\[
L_{n}^\kappa (t)= n!^{-1} e^t t^{-\kappa} d^n (e^{-t} t^{n +\kappa}) / dt^n ,\q t>0,
\]
and the functions 
\begin{equation}
u_{n}^\kappa (t)=\sqrt{\frac{n!}{\Gamma(n+1+\kappa)}} t^{\kappa /2}e^{-t/2}L_{n}^\kappa (t), \q n=0,1,\ldots,   
\label{eq:K3x}\end{equation}
form an orthonormal basis in the space $ L^2({\Bbb R}_{+})$. Therefore the operator  ${\bf U}_{\kappa}: l^2({\Bbb Z}_{+})\to  L^2({\Bbb R}_{+})$ defined by the formula
\begin{equation}
({\bf U}_{\kappa} g (t)=\sum_{n=0}^\infty g_{n} u_{n}^\kappa(t)
\label{eq:K3}\end{equation}
is unitary and hence
\[
({\bf U}_{\kappa}^{-1} f)_{n}=\int_{0}^\infty  f(t) u_{n}^\kappa(t) dt .
\]

Observe that if $g \in \ell_{0}$, then according to \e{eq:K3x} and  \e{eq:K3} we have
\begin{equation}
(\ov{ {\bf U}_0 g }\star {\bf U}_0 g )(t)= \sum_{n,m=0}^\infty  g_{m}\bar{g}_{n} e^{-t/2}\int_{0}^t L_{m}^0 (s) L_{n}^0 (t-s) ds.
\label{eq:KP1}\end{equation}
Putting together formulas (10.12.23) and (10.12.31) in \cite{BE}, we see that 
\begin{equation}
 \int_{0}^t L_{m}^0 (s) L_{n}^0 (t-s) ds= (n+m+1)^{-1} t L_{n+m}^1 (t).
\label{eq:KP1X}\end{equation}
Now to get relation \e{eq:KH} with ${\bf U} ={\bf U}_0$ and
\begin{equation}
 q_{n } =\frac{1}{n+1}\int_{0}^\infty  h  (t )   t L_{n}^1 (t ) e^{-t/2} dt,
\label{eq:K6}\end{equation}
we only have to multiply \e{eq:KP1} by $h(t)$ and integrate it in $t\in{\Bbb R}_{+}$. Since the operator  ${\bf U}_1$ is unitary, the last relation can formally be rewritten as
\begin{equation}
h  (t )=\sum_{n =0}^\infty q_{n } L_{n}^1 (t ) e^{-t/2}.
\label{eq:K5}\end{equation}

As usual, we consider $h(t)$ as a distribution.  The problem is that the functions $({\bf U}_0 g )(t)$ and  $t L_{n}^1 (t ) e^{-t/2} $ do not belong to the class $C_{0}^\infty ({\Bbb R}_{+})$, and hence the assumption $h\in C_{0}^\infty ({\Bbb R}_{+})'$ does not allow us to give a precise sense to formulas \e{eq:KH} and \e{eq:K6}. Therefore we introduce the set ${\cal X} \subset C^\infty ({\Bbb R}_{+})$ that consists of functions $\varphi(t)$ satisfying estimates
\begin{equation}
|\varphi^{(n)} (t)| \leq C_{n} t \q {\rm and}\q |\varphi^{(n)} (t)|\leq C_{n} e^{-\gamma t}
\label{eq:K55}\end{equation}
for some $\gamma< 1/2$ and all $n$. Since all functions $(\ov{ {\bf U}_0 g }\star {\bf U}_0 g )(t)$ and  $t L_{n}^1 (t ) e^{-t/2} $  belong to   ${\cal X}$, relations \e{eq:KP1} and \e{eq:KP1X} imply the following result.

\begin{theorem}\label{SIG2}
   Let $h\in{\cal X}'$, and let
   \[
   q_{n}={\pmb\la} h, \ti{u}^1_{n} {\pmb\ra} \q {\rm where}  \q \ti{u}^1_{n} (t)=   (n+1)^{-1}  t  L_{n}^1 (t ) e^{-t/2}.
   \]
Then for all elements $g \in {\ell}_{0}$   identity \e{eq:KH} holds with ${\bf U}  ={\bf U}_0 $.    
   \end{theorem}

Since the operator   ${\bf U}_{1}$ s unitary, it follows from \e{eq:K6} that
  \[
 \sum_{n =0}^\infty (n+1) |q_{n } |^2 = \int_{0}^\infty |h  (t )|^2 t dt.
\]
This relation simply means that the Hilbert-Schmidt norms of the operators $H$ and $Q$ related by formula \e{eq:K2}  are the same.

 We note that as shown in \cite{Y} (Theorem~3.8), the condition $h\in{\cal X}'$ is satisfied for all bounded operators.  For bounded  operators $H$ and $Q$,  identity \e{eq:KH} extends to all $g \in l^2({\Bbb Z}_{+})$ which yields  \e{eq:K2}.
Note that  $h\in{\cal X}'$ if 
\begin{equation}
\int_{0}^1 t |h(t)| dt+ \int_1^\infty e^{-\gamma t} |h(t)| dt <\infty
\label{eq:KP}\end{equation}
for some $\gamma <1/2$.  This assumption is by no means optimal although it even admits an exponential growth of $h(t)$ as $t\to\infty$. Even the condition $h\in L^1_{\rm loc} ({\Bbb R}_{+})$  is   not required. 

\begin{example}\label{SIex}
 Let  $h(t)=\d^{(k)} (t-t_{0})$ for some $k\in{\Bbb Z}_{+}$ and $t_{0}>0$;  then $h\in{\cal X}'$.  It follows from formula \e{eq:K6} that the   matrix elements of the corresponding Hankel operator $Q^{(k)}$ are given by the equality
\begin{equation}
 q_{n}^{(k)}=\frac{1}{n+1}(-1)^k  ( t L_{n}^1 (t ) e^{-t/2} )^{(k)}\big|_{t=t_{0}}.
\label{eq:KPQ}\end{equation}
 If $k=0$, then, as shown in \cite{Y}, the spectrum of the operator $H$ consists of the eigenvalues $-1,0,1$ of infinite multiplicity. According to Theorem~\ref{SIG2}, the spectrum of the  operator $Q^{(0)}$ is the same. Note that formula (10.15.1) in \cite{BE}  shows that
 \[
 L_{n}^1 (t_{0})= \pi^{-1/2} t_{0}^{-3/4} e^{t_{0}/2} n^{1/4}\cos (2\sqrt{n t_{0}}
-3\pi/4 )+O (n^{-1/4})
\]
so that even the boundedness of $Q^{(0)}$ does not look obvious. If $k\geq 1$, then (see \cite{Yf1}) the operators $H$ are unbounded and their spectra consist of eigenvalues accumulating both at $+\infty$ and $-\infty$. According to Theorem~\ref{SIG2}, the spectra of the Hankel operators $Q^{(k)}$ with matrix elements \e{eq:KPQ} possess the same properties.
   \end{example}

Recall that in {\it the class of bounded operators},  Hankel operators in $l^2 ({\Bbb Z}_{+})$ can be characterized by the commutation relation
\begin{equation}
QT=T^* Q
\label{eq:CR1}\end{equation}
where $T$ is the shift defined by the formula $(Tg)_{n}= g_{n-1}$ (with $g_{-1}=0$). Similarly, 
Hankel operators in $L^2 ({\Bbb R}_{+})$ can be characterized (see \cite{Y}, subs.~3.2, for details) by the commutation relation
\begin{equation}
H{\sf T}(\tau)={\sf T}(\tau)^* H, \q \forall \tau\geq 0,
\label{eq:CR2}\end{equation}
where $({\sf T}(\tau) f)(t)= f (t-\tau)$ for $t\geq \tau$ and $({\sf T}(\tau) f)(t)= 0$ for $t < \tau$.

  \begin{proposition}\label{CR}
  A bounded operator $H$ satisfies   \e{eq:CR2} if and only if the operator $Q={\bf U}_{0}^{-1} H {\bf U}_{0}$ satisfies  \e{eq:CR1}.
   \end{proposition}
   
     \begin{pf}
     As shown in \cite{Y} (see Corollary~3.5),  \e{eq:CR2} is equivalent to the relation
     $
H\Sigma= \Sigma^* H 
$
where
  \[
( \Sigma f)(t)=  e^{-t/2}\int_{0}^t e^{s/2} f(s) ds.
\]
Therefore we only have to verify that
$
{\bf U}_{0}  T  = ( I - \Sigma) {\bf U}_{0}. 
$
By definition  \e{eq:K3}, to that end we have to check the identity 
 \begin{equation}
u_{n+1}^0 (t)= u_{n }^0 (t) - e^{-t/2}\int_{0}^t e^{s/2} u_{n }^0 (s) ds 
\label{eq:CR6}\end{equation}
where $u_{n }^0 (t)$ are functions \e{eq:K3x}. Both sides of \e{eq:CR6} equal $1$ for $t=0$. The equality of their derivatives follows from the identity
\[
d \big( L_{n}^0 (t)    - L_{n+1}^0 (t) \big)/dt= L_{n}^0 (t)  
\]
(see formula (10.12.16) in \cite{BE}).
          \end{pf}

 It is possible to indicate the general form of unitary operators ${\bf U} : l^2({\Bbb Z}_{+})\to  L^2({\Bbb R}_{+})$ such that  operator  \e{eq:K2}   is a (bounded) Hankel operator  in $L^2({\Bbb R}_{+})$ if and only if $Q$ is a Hankel operator in $l^2({\Bbb Z}_{+})$. Let the dilation $D_{\rho}$, $\rho>0$,   be defined  in the space $L^2({\Bbb R}_{+})$ by the formula $( D_{\rho} f)(t) =\rho^{-1/2} f(\rho^{-1 }t)$, and let  the involution ${\cal J}$   be defined   in the space $L^2({\Bbb Z}_{+})$ by the formula $({\cal J}g)_{n} = (-1)^n g_{n}$. It is shown in the Appendix to
  \cite{Yf} that if an operator $\bf V$ is unitary in $L^2({\Bbb R}_{+})$ and ${\bf V}H {\bf V}^{-1}$ are Hankel operators  for all Hankel operators $H$, then necessarily either ${\bf V}=D_{\rho} $ or ${\bf V}=D_{\rho} {\bf U}_{0} {\cal J}  {\bf U}_{0} ^{-1}$ for some $\rho>0$.
  It follows   that  all unitary operators   ${\bf U} : l^2({\Bbb Z}_{+})\to  L^2({\Bbb R}_{+})$  possessing property \e{eq:K2}  admit one of the two following forms: ${\bf U}=D_{\rho}{\bf U}_{0}$ or ${\bf U}= D_{\rho}{\bf U}_{0} {\cal J} $ where $\rho>0$.
  
   Using this observation, we can  choose an arbitrary large $\gamma>0$ in definition \e{eq:K55} of the class $\cal X$. Then Theorem~\ref{SIG2}  remains true with    
   ${\bf U}=D_{\rho}{\bf U}_{0}$ for a suitable $\rho>0$. 
  
\medskip
 
 {\bf 7.2.}  
 Let us find a relation between matrix elements $q_{n}$ of a Hankel operator $Q$ and the  sigma-function $\sigma(\lambda)$ of the corresponding Hankel operator \e{eq:K2}.   Let us suppose that $\supp \sigma\subset [0, \infty)$ and substitute expression \e{eq:conv1} into formula \e{eq:K6}:
 \begin{equation}
 q_{n } = \frac{1}{n+1} \int_{0}^\infty d\lambda \sigma (\lambda) \int_{0}^\infty    t L_{n}^1 (t ) e^{-(1/2+\lambda)t} dt.
\label{eq:KL1}\end{equation}
According to formula (10.12.32) in \cite{BE} we have
 \[
  \int_{0}^\infty    t L_{n}^1 (t ) e^{-(1/2+\lambda)t} dt= \frac{n+1}{(\lambda+1/2)^2}\Big( \frac{ \lambda-1/2}{\lambda+1/2} \Big)^n ,
\]
and hence it follows from \e{eq:KL1}   that
 \[
 q_{n } =   \int_{0}^\infty \frac{\sigma (\lambda)}{(\lambda+1/2)^2}        \Big( \frac{ \lambda-1/2}{\lambda+1/2} \Big)^n d\lambda  .
\]
Introducing now the  function
 \begin{equation}
 \eta (\mu)= \sigma (\lambda)\q {\rm where} \q \mu=\frac{ \lambda-1/2}{\lambda+1/2} \in (-1,1),
\label{eq:KL3}\end{equation}
we obtain the representation \e{eq:KL5}.

  For the precise proof of   \e{eq:KL5}, we need only to justify the change of order of integrations in \e{eq:KL1}.
  It can be done by the Fubini theorem. We state only the simplest result which is however sufficient in many specific applications.
  
    \begin{proposition}\label{SIG1}
    Let the sigma-function $\sigma(\lambda)$ of a Hankel operator $H$ satisfy assumptions \e{eq:su}  and \e{eq:sig1}, and let  the function $\eta(\mu)$ be defined by formula \e{eq:KL3}. 
 Then   $Q={\bf U}_{0}^{-1} H {\bf U}_{0}$ is the Hankel operator
  in the space $l^2 ({\Bbb Z}_{+})$   with matrix elements  \e{eq:KL5}.    
   \end{proposition}
   
   Of course under the assumptions of this theorem  $|h(t)|Ê\leq C t^{-1}$, $\eta\in L^\infty (-1,1)$,  $q_{n}= O(n^{-1})$ and  the operators $H$ and $Q$ are bounded.
      
 \medskip
 
 {\bf 7.3.} 
 The method presented here gives, in principle, a constructive approach to the solution of the Hausdorff moment problem \e{eq:KL5}. We describe it in this subsection at a {\it formal} level.
 
  Given a sequence $q_{n}$, we first construct the kernel $h(t)$ by formula \e{eq:K5}. Then we find its sigma-function $\sigma(\lambda)$ by the inversion of the Laplace transform and, finally, we make the change of variables \e{eq:KL3}. The function $\eta(\mu)$ obtained (we also 
  call it the sigma-function of the Hankel operator $Q$) yields the solution of the   moment problem \e{eq:KL5}. In general,  $\eta $ is a distribution obtained from $\sigma\in {\cal Y}'$ by the change of variables \e{eq:KL3}, but $\eta (\mu)d\mu$   is a positive measure if $Q\geq 0$. We note (see Theorem~2.6.4 in \cite{AKH}) that original conditions for the solvability of the moment problem \e{eq:KL5} were formulated in rather different terms.
  
  Alternatively, we can exhibit an expression for the function $\eta(\mu)$, or rather for the Mellin transform of $\lambda^{-1/2}\sigma(\lambda)$, directly in terms of the coefficients $q_{n}$, avoiding the construction of the kernel $h(t)$:
  \begin{equation}
  \int_{0}^\infty \lambda^{-1 + i \xi} \sigma (\lambda) d \lambda= 2^{1- i\xi}  \sum_{n =0}^\infty i^{-n}  q_{n } P_{n}  (-\xi )  .
\label{eq:MP}\end{equation}
Here 
  \[
P_{n}(\xi)  =  i^n  \sum_{m=0}^n \frac{(-1)^m 2^m}{m !} C_{n+1}^{m+1} (1+i \xi)\cdots
(m + i \xi)
\]
 is the polynomial of degree $n$ known as the Meixner-Pollaczek polynomial; note that the term corresponding to $m=0$ in this sum is $1$ and $C_{n+1}^{m+1}$ are binomial coefficients.  Recall that $P_{n}(\xi) $, denoted also  $P^1_{n}(\xi; \pi/2) $ in \S 10.21 of \cite{BE},  are orthogonal   polynomials in the space $L^2 ({\Bbb R}; |\Gamma(1-i \xi) |^2 d\xi)$ related to the hypergeometric function $F(-n, 1-i\xi, 2;2)$ by formula (10.21.10).
  We give only a {\it formal} proof of relation \e{eq:MP}.  According to Lemma~\ref{MEan} for $\gamma=1$ we have 
  \[
\Gamma (1- i\xi) (M \Omega^{-1/2} \sigma)  (-\xi) = (M \Omega^{1/2} h)  (\xi) .
\]
Applying the operator $M \Omega^{1/2}$ to both sides of  equality \e{eq:K5}, we see that
 \[
(M \Omega h ) (\xi )=\sum_{n =0}^\infty  q_{n } \int_{0}^\infty t^{-i\xi} L_{n}^1 (t ) e^{-t/2} dt
\]
 where   
  \[
    \int_{0}^\infty t^{- i \xi} L^1_{n }(  t) e^{-t /2}dt= i^{-n} 2^{1- i\xi} \Gamma(1-i\xi) P_{n} (-\xi)
\]
according to  formula (10.12.33)  and expression (2.1.4)    for $F(-n, 1-i\xi, 2;2)$   in \cite{BE}.
Combining the formulas obtained, we  get relation \e{eq:MP}.

  \medskip
 
 {\bf 7.4.}  
 Let us  come back to quasi-Carleman operators $H$ with kernels   \e{eq:E1r} where we now suppose that $\alpha\geq 0$.  Our goal is to  calculate matrix elements $q_{n} = q_{n } (\alpha, k , r) $ of Hankel operators $Q={\bf U}_{0}^{-1}H{\bf U}_{0}$. Let us proceed   from formula   \e{eq:KL5} for $q_{n}$ in terms of the sigma-function. We suppose that $k\not\in{\Bbb Z}_{+}$ since for $k \in{\Bbb Z}_{+}$ the operators $Q$ have finite rank and the well-known formula for $q_{n}$ is, for example, a direct consequence of \e{eq:dii}.  Putting together relations 
 \e{eq:bbr7}  and  \e{eq:KL3}, we obtain the following result.
 
   \begin{proposition}\label{SIM}
  Let $k\not\in{\Bbb Z}_{+}$. In the case $r=0$ assume additionally that   $k> -2$. Then the  matrix elements 
   $q_{n}$ of the operator $Q={\bf U}_{0}^{-1}H{\bf U}_{0}$ are determined by relation  \e{eq:KL5} where
  \begin{equation}
\eta(\mu)= (\alpha+1/2)^{-k-1} \Gamma (-k)^{-1} (1-\mu)^{k+1}(\mu-\gamma)_{+}^{-k-1}\exp\big( \! -r (\alpha+1/2)\frac{\mu -\gamma} {1-\mu}\big)
\label{eq:KS}\end{equation}
 with
  \[
 \gamma=\frac{ \alpha-1/2}{\alpha+1/2}\in [-1,1).
\]
     \end{proposition}
     
     The case $r=0$ is particularly simple. If $\alpha=1/2$, then 
 it follows from  \e{eq:KL5}, \e{eq:KS} that   
\[
 q_{n } = \frac{\Gamma (k+2) \Gamma (n-k) }{  \Gamma (-k)  \Gamma (n+2)}.  
\]
 If $k=-1$, then 
 \begin{equation}
q_{n}= \int_{\gamma}^1 \mu^n d\mu= \frac{1-\gamma^{n+1}}{n+1}.
\label{eq:KMzr}\end{equation}
If $\alpha=0$, then $\gamma=-1$ and we recover of course the matrix elements  of the Carleman operator. If $\alpha>0$, then $\gamma \in (-1,1)$ and we obtain the generalized (but different from those considered by M.~Rosenblum in \cite{Ro}) Hilbert matrices. They reduce to the standard Hilbert matrix for $\gamma=0$.

Alternatively, we could proceed from  
 formula \e{eq:K6} for $q_{n}$ in terms of the kernel $h(t)$. In the particular case $r=0$, using formula
 (10.12.33) in \cite{BE},  we can express the coefficients $q_{n }  $  via the hypergeometric function $F$:
 \begin{equation}
 q_{n } = q_{n} (\beta, k)=  \Gamma (2+k) \beta^{2+k}  F (-n, 2+k, 2; \beta)  ,\q \beta= (\alpha+1/2)^{-1}\in (0,2] .
\label{eq:KM}\end{equation}
 
 Observe that if $r=0$, then  formulas \e{eq:KL5} as well as  \e{eq:K6}  make sense for $k> -2$ only while  the Hankel quadratic form ${\pmb\la }h,\bar{f}\star f {\pmb\ra }$  is well defined on $f\in C_{0}^\infty ({\Bbb R}_{+})$ for all  kernels $h(t)= t^k e^{-\alpha t}$. 
Thus the example of quasi-Carleman operators  shows that   considerations of Hankel operators
   in the spaces $L^2 ({\Bbb R}_{+})$ and $l^2 ({\Bbb Z}_{+})$ are not always equivalent.
   
Note that, for all $\rho>0$, Hankel operators $H$ and $H_{\rho}$ with kernels $h(t)$ and $h_{\rho}(t)= \rho h(\rho t)$ are unitarily equivalent (by the dilation transformation). In particular, all Hankel operators $H_{\rho}$  with kernels   $h_{\rho}(t)=\rho^{1+k}t^k e^{-\rho t/2}$ are unitarily equivalent to each other for all $\rho>0$. This implies the following assertion.

  \begin{proposition}\label{QCM}
  Let $Q(\beta, k)$  be  the
Hankel operator in the  space $l^2 ({\Bbb Z}_{+})$ with the  matrix elements  \e{eq:KM} where $k>-2$. Then the operators 
 $(-1+2/\beta)^{1+k} Q( \beta, k)$  are unitarily equivalent to each other for all $\beta \in (0,2)$. In particular $($for $k=-1)$, the  generalized  Hilbert matrices determined by formula \e{eq:KMzr} are unitarily equivalent to each other for all $\gamma \in (-1,1)$. Thus their spectra are absolutely continuous, simple and coincide with the interval $[0,\pi]$.
     \end{proposition}

This result does not look obvious in the discrete representation $l^2 ({\Bbb Z}_{+})$, but it becomes quite transparent  after the transformation of the problem    into the space $L^2 ({\Bbb R}_{+})$.

 \medskip
 
 {\bf 7.5.} 
 Our next goal is to find the asymptotics of matrix elements $q_{n}  $ of the quasi-Carleman operators as $n\to\infty$. It  easily follows from formula \e{eq:KS} that for any $a\in (0,1)$
   \[
 \int_{-a}^a \eta(\mu) \mu^n d\mu= O (b^n), \q \forall b>a,
\]
 so that the asymptotics of   $q_{n}  $ is determined by neighborhoods of the points $\mu=\pm 1$ in the integral representation \e{eq:KL5}.
 
 Consider first the point $\mu=-1$. If $\alpha>0$, that is $\gamma>-1$, then function \e{eq:KS} equals zero in a neighborhood of the point $-1$. So the contribution of this point to  the asymptotics of   $q_{n}  $ is  also zero. If $\alpha=0$, that is $\gamma= -1$, then it follows from formula \e{eq:KS} that
  \[
\eta(\mu)= 4^{k+1} \Gamma (-k)^{-1}  (\mu+1)^{-k-1} \big( 1+ O (\mu+1)\big)
\]
as $\mu\to -1$. So we have
  \begin{multline}
q_{n}^{(-)}:=  \int_{-1}^0 \eta (\mu ) \mu^n d\mu
\\
= 4^{k+1} \Gamma (-k)^{-1} \int_{-1}^0 (\mu+1)^{-k-1}\mu^n d\mu+ O\big(\int_{-1}^0 (\mu+1)^{-k}|\mu|^n d\mu\big ).
 \label{eq:KS2}\end{multline}
 Note that
 \[
\int_{-1}^0 (\mu+1)^{-k-1}\mu^n d\mu= (-1)^n \frac{ \Gamma (-k) \Gamma (n+1)}{    \Gamma (n-k+1)} = (-1)^n 
\Gamma (-k) n^k (1+O (n^{-1})) 
\] 
where we have used the asymptotic formula   (1.18.4) in  \cite{BE} for the ratio of the gamma functions. Thus according to
\e{eq:KS2} we have
\begin{equation}
q_{n}^{(-)}= (-1)^n  4^{k+1} n^k (1+O (n^{-1})). 
\label{eq:KS3}\end{equation}

 Next, we consider a neighborhood of the point $\mu= 1$. If $r>0$, then function  \e{eq:KS} exponentially tends to zero as $\mu\to 1$ so that the contribution $q_{n}^{(+)}$ of this point to $q_{n}$ is negligible.  If $r=0$, then it follows from formula \e{eq:KS} that
  \begin{equation}
\eta(\mu)=   \Gamma (-k)^{-1}  (1-\mu)^{k+ 1} \big( 1+ O (1-\mu)\big)
\label{eq:KS1+}\end{equation}
as $\mu\to 1$. So we have
  \[
q_{n}^{(+)}:=  \int_0^1 \eta (\mu ) \mu^n d\mu= \Gamma (-k)^{-1}  \int_{0}^1   (1-\mu)^{k+ 1} \mu^n d\mu
+O \big( \int_{0}^1   (1-\mu)^{k+ 2} \mu^n d\mu\big) 
 \]
 Calculating again the integrals here in terms of the beta function and using   formula   (1.18.4) in  \cite{BE}, we find that
   \begin{equation}
q_{n}^{(+)} =  \Gamma (-k)^{-1} \Gamma (k+2)  n^{-k-2} \big( 1+ O (n^{-1})\big).
\label{eq:KS2+}\end{equation} 

Let us put together the results obtained.

   \begin{proposition}\label{SIMA}
  Let the assumptions of Proposition~\ref{SIM} hold. 
  If $\alpha>0$ and $r>0$, then  
 the  matrix elements 
   $q_{n}$ of the Hankel operator $Q $ decay faster than any power of $n^{-1}$ as $n\to \infty$.
   If $\alpha=0$ but  $r>0$, then  the asymptotics of   $q_{n}$ 
 is given by formula \e{eq:KS3} where $q_{n}= q_{n}^{(-)}$. If $r=0$ but  $\alpha > 0$, then  the asymptotics of   $q_{n}$ 
 is given by formula \e{eq:KS2+} where $q_{n}= q_{n}^{(+)}$. Finally, if $\alpha=r =0$, then 
    \begin{equation}
q_{n}  = (-1)^n  4^{k+1} n^k (1+O (n^{-1})) + \Gamma (-k)^{-1} \Gamma (k+2)  n^{-k-2} \big( 1+ O (n^{-1})\big).
\label{eq:KSL}\end{equation}
       \end{proposition}

   \begin{remark}\label{SIMAx}
Of course if $k\in (-2, -1)$ (if $k>-1$), then the first  (the second) term in \e{eq:KSL} can be neglected. If $k=-1$, then both terms in \e{eq:KSL}  have the same order.
       \end{remark}
  
We emphasize that under the assumptions of Proposition~\ref{SIMA} the sequence $q_{n}$ does not necessarily tend to $0$ and the operator $Q$ may be unbounded.

\appendix{}

\section{A generalization of the Bernstein theorem}
 
  Our proof of Theorem~\ref{Bern} will be divided in a  series of simple lemmas. For  an arbitrary $\varphi\in C_{0}^\infty({\Bbb R}_{+})$, we set
 \begin{equation}
  \eta =\varphi\star \bar{\varphi}
 \label{eq:eta}\end{equation}
  and define the distribution
  \begin{equation}
h_{\varphi}(t)=\int_{0}^\infty h(s) \eta(s-t) ds={\pmb \la}h, \eta(\cdot-t) {\pmb \ra}
\label{eq:Bern2}\end{equation}
which is actually    a continuous function of $t>0$.
It follows from \e{eq:eta} and \e{eq:Bern2} that
  \begin{equation}
{\pmb \la}h_\varphi , f {\pmb \ra}= {\pmb \la}h  , \varphi\star \bar{\varphi} \star f {\pmb \ra}
\label{eq:Bern2a}\end{equation}
 for an arbitrary $f\in C_{0}^\infty({\Bbb R}_{+})$. 

Let us   check that
 \begin{equation}
\int_{0}^\infty\int_{0}^\infty h_\varphi(\tau+\sigma) g(\sigma)\overline{g(\tau)} d\tau d\sigma\geq 0
\label{eq:Bern3}\end{equation}
for all $g\in C_{0}^\infty({\Bbb R}_{+})$. Using \e{eq:eta} and  \e{eq:Bern2}, we can rewrite the last  integral   as 
\[
\int_{0}^\infty\int_{0}^\infty \int_{0}^\infty\int_{0}^\infty 
h (t+s+\tau+\sigma)  \varphi (s) \overline{\varphi (t)} g(\sigma)\overline{g(\tau)} dt ds d\tau d\sigma.
\]
Making here the changes of variables $x=t+\tau$, $y=s+\sigma$, we see that this expression equals
 \begin{equation}
\int_{0}^\infty\int_{0}^\infty  
h (x+y)  \psi(y) \overline{\psi(x)} dx dy
\label{eq:Bern5}\end{equation}
where 
\[
\psi(x)= \int_{0}^x g(x-t)\varphi (t)dt.
\]
Since $\psi\in C_{0}^\infty({\Bbb R}_{+})$, expression \e{eq:Bern5} is positive by the condition \e{eq:BernH}. This proves \e{eq:Bern3}.

Thus applying the Bernstein theorem on exponentially convex functions to the function $h_\varphi(t)$, we obtain the following intermediary  result.

\begin{lemma}\label{Bern1}
For an arbitrary  $\varphi \in C_{0}^\infty({\Bbb R}_{+})$, let the function $h_\varphi(t)$ be defined by formulas \e{eq:eta}, \e{eq:Bern2}. Then under the assumption  \e{eq:BernH}, there exists a positive measure $M_\varphi$ on $\Bbb R$ such that
\begin{equation}
h_\varphi(t)=\int_{-\infty}^\infty e^{-t\lambda} dM_\varphi(\lambda)  
\label{eq:Bern4}\end{equation}
where the integral is convergent for all $t>0$.
 \end{lemma}
 
 Let us now compare the measures $M_\varphi$ corresponding to different functions $\varphi$.

\begin{lemma}\label{Bern2}
For arbitrary functions $\varphi_1, \varphi_2\in C_{0}^\infty({\Bbb R}_{+})$,
   $j=1,2$, and all
$\lambda\in{\Bbb R}$ the relation  
\begin{equation}
|({\sf L}\varphi_{2})(\lambda)|^2 d M_{\varphi_{1}}(\lambda)
=|({\sf L}\varphi_{1})(\lambda)|^2 d M_{\varphi_{2}}(\lambda)  
\label{eq:Bern6}\end{equation}
holds.
 \end{lemma}

\begin{pf}
Let us proceed from definition \e{eq:Bern2a} which,   for an arbitrary $f\in C_{0}^\infty ({\Bbb R}_{+})$, yields   relations
 \[
{\pmb \la}h_{\varphi_{1}} , \eta_{2}\star f {\pmb \ra} = {\pmb \la}h  , \eta_1\star (\eta_{2}\star f) {\pmb \ra},
\q
{\pmb \la}h_{\varphi_2} , \eta_1\star f {\pmb \ra} = {\pmb \la}h  , \eta_{2}\star (\eta_1\star f) {\pmb \ra}
 \]
where notation \e{eq:eta} has been used.
 Since
$
   \eta_1\star (\eta_{2}\star f)=\eta_{2}\star (\eta_1\star f),
  $
   it follows that
    \begin{equation}
{\pmb \la}h_{\varphi_{1}} , \eta_{2}\star f {\pmb \ra} =
{\pmb \la}h_{\varphi_2} , \eta_1\star f {\pmb \ra} .
  \label{eq:Bern2c}\end{equation}
  
According to Lemma~\ref{Bern1}  we have
    \begin{multline*}
{\pmb \la}h_{\varphi_{1}} , \eta_{2}\star f {\pmb \ra} =
\int_{-\infty}^\infty \ov{({\sf L}( \eta_{2}\star f ))(\lambda)} d M_{\varphi_{1}} (\lambda)
=\int_{-\infty}^\infty \ov{({\sf L} f)(\lambda)   ({\sf L} \eta_{2})(\lambda)   }d M_{\varphi_{1}} (\lambda)
 \\
=
\int_0^\infty dt \ov{ f(t) }\int_{-\infty}^\infty e^{-t\lambda}  \ov{ ({\sf L} \eta_{2})(\lambda) } d M_{\varphi_{1}} (\lambda).
 \end{multline*}
   The exactly similar identity is true for ${\pmb \la}h_{\varphi_2} , \eta_1\star f {\pmb \ra}$.
  Hence it follows from equality \e{eq:Bern2c} that
 \[ 
\int_0^\infty dt \ov{ f(t) }\int_{-\infty}^\infty  e^{-t\lambda}  \ov{ ({\sf L} \eta_{2})(\lambda) } d M_{\varphi_{1}} (\lambda)
=\int_0^\infty dt \ov{f(t)} \int_{-\infty}^\infty e^{-t\lambda} \ov{  ({\sf L} \eta_1)(\lambda) } d M_{\varphi_2} (\lambda).
 \]
 This ensures that for all $t>0$
 \[
  \int_{-\infty}^\infty e^{-t\lambda}  \ov{ ({\sf L} \eta_{2})(\lambda) } d M_{\varphi_{1}} (\lambda)=  \int_{-\infty}^\infty e^{-t\lambda}  \ov{ ({\sf L} \eta_1)(\lambda) } d M_{\varphi_2} (\lambda)
\]
because $f\in C_{0}^\infty ({\Bbb R}_{+})$ is arbitrary. 
Therefore by the uniqueness theorem for Laplace integrals we have
 \[
   \ov{ ({\sf L} \eta_{2})(\lambda)  }d M_{\varphi_{1}} (\lambda)= \ov{   ({\sf L} \eta_1)(\lambda) } d M_{\varphi_2} (\lambda), \q \forall \lambda\in{\Bbb R}.
\]
Since 
$
 ({\sf L} \eta_{j})(\lambda) = |({\sf L} \varphi_{j})(\lambda)|^2,\q j=1,2,
$
 this is equivalent to   identity \e{eq:Bern6}.  
\end{pf}

Now we define the measure $dM(\lambda)$ on $\Bbb R$  by the relation
   \begin{equation}
 dM(\lambda)=  |({\sf L} \varphi )(\lambda)|^{-2} dM_{\varphi}(\lambda) .
  \label{eq:Ber}\end{equation}
In view of  Lemma~\ref{Bern2} this definition does not depend on the choice of the function $\varphi\in C_{0}^\infty ({\Bbb R}_{+})$.

\begin{lemma}\label{Bern11}
 The measure \e{eq:Ber} satisfies estimates \e{eq:Bern1x}.
 \end{lemma}

\begin{pf}
We proceed from estimates \e{eq:Bern1x} on the measures   
 $dM_{\varphi}(\lambda)$. Observe that if  $\varphi\in C_{0}^\infty ({\Bbb R}_{+})$,  $\varphi\neq 0$, $\varphi(t)\geq 0$ for $t < t_{0}/2$ and  $\varphi(t)= 0$ for $t\geq t_{0}/2$, then 
 \begin{equation}
 ({\sf L} \varphi )(\lambda)\geq c(\varphi) e^{-t_{0}\lambda/2},
 \q c(\varphi)=\int_{0}^\infty \varphi(t)dt>0.
  \label{eq:BerBer}\end{equation}
 It follows from \e{eq:Ber} that 
 \[
 \int_{0}^\infty e^{-t\lambda} dM(\lambda)\leq  c(\varphi)^{-2}
  \int_{0}^\infty e^{-(t -t_{0})\lambda} dM_{\varphi}(\lambda).
  \]
  For an arbitrary $t>0$, we can choose $t_{0}$ so small that $t-t_{0}>0$, and hence the  integral on the right is convergent. This yields the first estimate \e{eq:Bern1x} for the measure $dM(\lambda)$.
  
  The second  estimate \e{eq:Bern1x} is even simpler because  if  $\varphi\in C_{0}^\infty ({\Bbb R}_{+})$,  $\varphi\neq 0$ and $\varphi(t)\geq 0$,   then $ ({\sf L} \varphi )(-\lambda)\geq c(\varphi) $ with $c(\varphi)$ defined in \e{eq:BerBer} for all $\lambda\geq 0$.
 It follows again from \e{eq:Ber}  that 
   \[
   \int_{0}^\infty e^{t\lambda} dM(-\lambda)\leq c(\varphi)^{-2} 
   \int_{0}^\infty e^{t\lambda} dM_{\varphi}(-\lambda)  
\]
where the  integral on the right is convergent for an arbitrary large $t>0$.
\end{pf}

 It remains to verify representation  \e{eq:Bern1}, or equivalently \e{eq:Bern1X},  for the measure $dM(\lambda)$  defined by  relation   \e{eq:Ber}. 
 
\begin{lemma}\label{Bern12}
Representation    \e{eq:Bern1X} is true for functions
 $F =\eta\star f $
where $\eta$ is function \e{eq:eta} and $f\in C_{0}^\infty({\Bbb R}_{+})$ is arbitrary.
 \end{lemma}

\begin{pf} 
  It follows from relations \e{eq:Bern2a},    \e{eq:Bern4} and definition  \e{eq:Ber} of the measure $dM (\lambda)$ that
  \begin{equation}
{\pmb \la}  h, F {\pmb \ra}={\pmb \la}  h_{\varphi}, f {\pmb \ra}=\int_{-\infty}^\infty \ov{({\sf L} f)(\lambda) } dM_{\varphi}(\lambda)= \int_{-\infty}^\infty \ov{({\sf L} f)(\lambda)} | ({\sf L} \varphi)(\lambda)|^2 dM (\lambda) .
\label{eq:fgp}\end{equation}
According to  property \e{eq:gg1} of the Laplace transform we have
 \[
  ({\sf L} F)(\lambda) = ({\sf L} f)(\lambda) | ({\sf L} \varphi)(\lambda)|^2.
\]
Therefore equality \e{eq:fgp} yields \e{eq:Bern1X}.
 \end{pf} 
 
 Finally, we extend representation  \e{eq:Bern1X} to all functions $F \in C_{0}^\infty({\Bbb R}_{+})$. Let us choose $\omega =\bar{\omega}\in C_{0}^\infty({\Bbb R}_{+})$ such that
 $
 \int_{0}^\infty \omega(t)dt =1
 $
 and set
 \begin{equation}
 \varphi_{\varepsilon}(t)= \varepsilon^{-1} \omega(\varepsilon^{-1}  t).
\label{eq:fgps}\end{equation}
It follows from Lemma~\ref{Bern12}  that
 \begin{equation}
{\pmb \la}  h,   \varphi_{\varepsilon} \star  {\varphi}_{\varepsilon} \star F {\pmb \ra}= 
\int_{-\infty}^\infty \ov{ ({\sf L} F)(\lambda)} | ({\sf L} \varphi_{\varepsilon})(\lambda)|^2 dM (\lambda).
\label{eq:fgp1}\end{equation}

Let us pass here to the limit $\varepsilon\to 0$. 

\begin{lemma}\label{Bern13}
Let $F \in C_{0}^\infty({\Bbb R}_{+})$. Then
 \begin{equation}
 \varphi_{\varepsilon} \star \varphi_{\varepsilon} \star F \to F
\label{eq:pas}\end{equation}
as $\varepsilon\to 0$ in the space $C_{0}^\infty({\Bbb R}_{+})$.
 \end{lemma}

\begin{pf}
 Since
  \[
  (\varphi_{\varepsilon} \star  {\varphi}_{\varepsilon})(t)= \varepsilon^{-1} \zeta(\varepsilon^{-1}  t)
  \]
  where $\zeta=\omega\star \omega$, we have
  \begin{align*}
  ( \varphi_{\varepsilon} \star \varphi_{\varepsilon} \star F)(t)&=
  \varepsilon^{-1} \int_{0}^t \zeta(\varepsilon^{-1}  s) F(t-s) ds
\\  =
  \int_{0}^{t/\varepsilon}& \zeta(\sigma) F (t-\varepsilon\sigma) d\sigma\to F(t)
   \int_{0}^{\infty} \zeta(\sigma)   d\sigma = F(t)
\end{align*}
as $\varepsilon\to 0$. Similar relations are of course also true for all derivatives in $t$.
Since the supports of $\varphi_{\varepsilon}$ are  small, the supports of all functions
$\varphi_{\varepsilon} \star \varphi_{\varepsilon} \star F$ are contained in a common interval $[t_{1},t_{2}]\in{\Bbb R}_{+}$. This leads to \e{eq:pas}.
\end{pf}

Thus  the left-hand side of \e{eq:fgp1} converges to the left-hand side of  \e{eq:Bern1X}.

\begin{lemma}\label{Bern14}
Let $F \in C_{0}^\infty({\Bbb R}_{+})$. 
Then the right-hand side of \e{eq:fgp1} converges to the right-hand side of  \e{eq:Bern1X}.
 \end{lemma}
 
\begin{pf}
It follows from \e{eq:fgps} that
 \begin{equation}
 ( {\sf L}\varphi_{\varepsilon} )(\lambda) = \int_{0}^\infty e^{-\varepsilon \lambda s} \omega(s)ds,
\label{eq:pas11}\end{equation}
and hence
 $ ( {\sf L}\varphi_{\varepsilon} )(\lambda) \to 1$ as $\varepsilon\to 0$
  for all $\lambda\in {\Bbb R}$. Moreover, if $\supp \omega\in [\omega_{1}, \omega_{2}]$, then  function \e{eq:pas11} is bounded by $Ce^{-\varepsilon\lambda \omega_{1}}$ for $\lambda\geq 0$ and by $Ce^{\varepsilon |\lambda| \omega_{2}}$ for $\lambda\leq 0$. Recall also that the measure $dM(\lambda)$ satisfies estimates \e{eq:Bern1}. Thus by the dominated convergence theorem, the right-hand side of \e{eq:fgp1} converges as $\varepsilon\to 0$ to  the right-hand side of  \e{eq:Bern1X}
  \end{pf}

Putting together relation \e{eq:fgp1} with Lemmas~\ref{Bern13} and \ref{Bern14}, 
we conclude the proof of Theorem~\ref{Bern}.


 \end{document}